\documentclass[fleqn,12pt,a4paper]{article}

\usepackage{ddwmacs,latexsym,amssymb}
\usepackage[dvips,xdvi]{graphics,color}
\usepackage{epsfig}

\title{%
  On the Links--Gould Invariant of Links
}

\date{}

\author{%
  David~~De Wit%
  \footnote{
    David De Wit and Jon R Links:
    Department of Mathematics, The University of Queensland.
    Q, 4072, Australia.
    email: \texttt{ddw@maths.uq.edu.au}, \texttt{jrl@maths.uq.edu.au}.
  },
  Louis H Kauffman%
  \footnote{
    Louis H Kauffman:
    Department of Mathematics, Statistics and Computer Science,
    The University of Illinois at Chicago.
    851 South Morgan Street, Chicago IL, 60607-7045, USA.
    email: \texttt{kauffman@uic.edu}.
  }~
  and
  Jon R Links$^*$
}

\begin{document}

\maketitle

\begin{center}
  \begin{picture}(0,0)%
\epsfig{file=BuddhistKnot.pstex}%
\end{picture}%
\setlength{\unitlength}{0.00083300in}%
\begingroup\makeatletter\ifx\SetFigFont\undefined
\def\x#1#2#3#4#5#6#7\relax{\def\x{#1#2#3#4#5#6}}%
\expandafter\x\fmtname xxxxxx\relax \def\y{splain}%
\ifx\x\y   
\gdef\SetFigFont#1#2#3{%
  \ifnum #1<17\tiny\else \ifnum #1<20\small\else
  \ifnum #1<24\normalsize\else \ifnum #1<29\large\else
  \ifnum #1<34\Large\else \ifnum #1<41\LARGE\else
     \huge\fi\fi\fi\fi\fi\fi
  \csname #3\endcsname}%
\else
\gdef\SetFigFont#1#2#3{\begingroup
  \count@#1\relax \ifnum 25<\count@\count@25\fi
  \def\x{\endgroup\@setsize\SetFigFont{#2pt}}%
  \expandafter\x
    \csname \romannumeral\the\count@ pt\expandafter\endcsname
    \csname @\romannumeral\the\count@ pt\endcsname
  \csname #3\endcsname}%
\fi
\fi\endgroup
\begin{picture}(2423,1928)(2130,-1625)
\end{picture}

\end{center}

\begin{abstract}
  \noindent
  We introduce and study in detail an invariant of $(1,1)$ tangles.
  This invariant, derived from a family of four dimensional
  representations of the quantum superalgebra $U_q[gl(2|1)]$, will be
  referred to as the Links--Gould invariant. We find that our invariant
  is distinct from the Jones, HOMFLY and Kauffman polynomials
  (detecting chirality of some links where these invariants fail), and
  that it does not distinguish mutants or inverses. The method of
  evaluation is based on an abstract tensor state model for the
  invariant that is quite useful for computation as well as theoretical
  exploration.
\end{abstract}


\section{Introduction}

Since the discovery of the Jones polynomial \cite{Jones:85}, several
new invariants of knots, links and tangles have become available due to
the development of sophisticated mathematical techniques. Among these,
the quantum algebras as defined by Drinfeld \cite{Drinfeld:87} and
Jimbo \cite{Jimbo:85}, being examples of quasi-triangular Hopf
algebras, provide a systematic means of solving the Yang--Baxter
equation and in turn may be employed to construct representations of
the braid group. From each of these representations, a prescription
exists to compute invariants of oriented knots and links
\cite{Reshetikhin:87,Turaev:88,ZhangGouldBracken:91b}, from which the
Jones polynomial is recoverable using the simplest quantum algebra
$U_q\[sl(2)\]$ in its minimal (2-dimensional) representation.

\vfill

\pagebreak

From such a large class of available invariants, it is natural to ask
if generalisations exist, with the view to gaining a classification.
One possibility is to look to multiparametric extensions in order
to see which invariants occur as special cases. A notable example is
the HOMFLY invariant \cite{FreydYetterHosteLickorishMilletOcneanu:85}
which includes both the Jones and Alexander--Conway invariants
\cite{Alexander:23,Conway:70} as particular cases as well as the
invariants arising from minimal representations of $U_q\[sl(n)\]$
\cite{Turaev:88}. Another is the Kauffman polynomial which includes the
Jones invariant as well as those obtained from the quantum algebras
$U_q\[o(n)\]$ and $U_q\[sp(2n)\]$ in the $q$-deformations of the
defining representations \cite{Turaev:88}.

The work of Turaev and Reshetikhin \cite{ReshetikhinTuraev:90} shows
that the algebraic properties of quantum algebras are such that an
extension of this method to produce invariants of oriented tangles is
permissible. A tangle diagram is a link diagram with free ends.  An
associated invariant takes the form of a tensor operator acting on a
product of vector spaces. Zhang \cite{Zhang:95} has extended this
formalism to the case of quantum superalgebras which are
$\mathbb{Z}_2$-graded generalisations of quantum algebras.

Since quantum superalgebras give rise to nontrivial one-parameter
families of irreducible representations, it is possible to utilise them
for the construction of two variable invariants. This was first shown
by Links and Gould \cite{LinksGould:92b} for the simplest case using the
family of four dimensional representations of $U_q\[gl(2|1)\]$. It was
also made known that a one variable reduction of this invariant
coincides with a one variable reduction of the Kauffman polynomial by
the use of the Birman--Wenzl--Murakami algebra. Extensions to more
general representations of quantum superalgebras are discussed in
\cite{GouldLinksZhang:96b}.

Thus far, little has been investigated with regard to the Links--Gould
invariant. Here we report on some properties and behaviour. The method
of evaluating the invariant involves a prior construction of the
quantum $R$-matrix associated with a family of four dimensional
representations.  Having obtained this matrix, the construction of the
invariant follows from properties of ribbon Hopf (super)algebras and
their representations. Here we consider the invariants of $(1,1)$
tangles for the following reason: for invariants derived from
representations of quantum superalgebras with zero $q$-superdimension,
the corresponding invariant is also zero.  If the representation is
irreducible, the quantum superalgebra symmetry of the procedure ensures
that the invariant of $(1,1)$ tangles takes the form of some scalar
multiple of the identity matrix.  (See \cite{ReshetikhinTuraev:90} for
a discussion of this symmetry.) We take this scalar to be the
invariant.

In this paper, we prove that the Links--Gould invariant is not able to
distinguish between mutant links (\secref{KTMutants}), nor is it able
to distinguish a knot from its inverse (i.e. from the knot obtained by
reversing the orientation) (\propref{inversion}). However it is good at
distinguishing some knots and links from their mirror images (see
Propositions \ref{prp:chiralmeanspalindromic} and
\ref{prp:Rhadinocentrus}), and it is distinguished from the HOMFLY and
Jones polynomials by this behaviour (see
\secref{NineFortyTwoandTenFortyEight} for specific examples). We
provide many examples and a complete description of the state model for
the invariant in abstract tensor form. This description of the
invariant directly facilitates the construction of a computer program
in \textsc{Mathematica} for calculation of the invariant.



\section{Construction of the $ R $ Matrix}
\seclabel{Rmatrixdef}

We consider the family of four dimensional representations of the
quantum superalgebra $ U_q \[ gl \( 2 | 1 \) \] $, which depend on a
free complex parameter $ \alpha $.
This superalgebra has $ 7 $ simple generators
$
  \{
    {E^1}_1, {E^2}_2, {E^3}_3,
    {E^1}_2, {E^2}_1,
    {E^2}_3, {E^3}_2
  \}
$
on which we define a $ \BZ_2 $ grading in
terms of the natural grading on the indices $ \[ 1 \] = \[ 2 \] = 0 $,
$ \[ 3 \] = 1 $ by:
\be
  \[ {E^i}_j \] = \[ i \] + \[ j \]
  \qquad \( \mod~2 \).
\ee

The $ U_q \[ gl \( 2 | 1 \) \] $ generators satisfy the commutation
relations:
\be
  \[ {E^1}_2, {E^2}_1 \]
  \eq
  {\[ {E^1}_1 - {E^2}_2 \]}_q
  \\
  \{ {E^2}_3, {E^3}_2 \}
  \eq
  {\[ {E^2}_2 + {E^3}_3 \]}_q
  \\
  \[ {E^i}_i, {E^j}_k \]
  \eq
  {\delta^j}_i {E^i}_k
  -
  {\delta^i}_k {E^j}_i,
  \qquad
  i, j, k = 1, 2, 3,
\ee
where $ \[~,~\] $ and $ \{~,~\} $ denote the usual commutator and
anticommutator, respectively and
we have employed the \emph{$ q $ bracket}, defined for a wide
class of objects $ x $ by:
\be
  {[x]}_q
  \defeq
  \frac{q^x - q^{-x}}{q - q^{-1}}.
\ee

Let $ \{ \ket{i} \}_{i = 1}^4 $ denote a basis for
the four dimensional $ U_q \[ gl \( 2 | 1 \) \] $ module
$ V $.  Consistent with the $ \mathbb{Z}_2 $ grading
on $ U_q \[ gl \( 2 | 1 \) \] $, we grade the basis states by:
\be
  \[ \ket{1} \] = \[ \ket{4} \] = 0,
  \qquad
  \[ \ket{2} \] = \[ \ket{3} \] = 1.
\ee
We define a dual basis $ \{ \bra{i} \}_{i = 1}^4 $; in component form,
these are represented by the transpose complex conjugates
of the original basis:
$ \bra{i} = \overline{\ket{i}}^t \equiv \ket{i}^\dagger $. Then:
$ \bra{i} \ket{j} \equiv \langle i | j \rangle = \delta_{i j} $.
In terms of these dual bases, we define a representation $ \pi $ of the
$ U_q \[ gl \( 2 | 1 \) \] $ generators; their action on the basis
vectors is given by:
\be
  \pi \( {E^1}_1 \)
  \eq
  - \ket{2} \bra{2} - \ket{4} \bra{4}
  \\
  \pi \( {E^2}_2 \)
  \eq
  - \ket{3} \bra{3} - \ket{4} \bra{4}
  \\
  \pi \( {E^3}_3 \)
  \eq
  \alpha \ket{1} \bra{1}
  +
  (\alpha+1) \( \ket{2} \bra{2} + \ket{3} \bra{3} \)
  +
  (\alpha+2) \ket{4} \bra{4}
  \\
  \pi \( {E^1}_2 \)
  \eq
  -
  \ket{3} \bra{2}
  \\
  \pi \( {E^2}_1 \)
  \eq
  -
  \ket{2} \bra{3}
  \\
  \pi \( {E^2}_3 \)
  \eq
  {[\alpha]}_q^{1/2} \ket{1} \bra{3}
  -
  {[\alpha+1]}_q^{1/2} \ket{2} \bra{4}
  \\
  \pi \( {E^3}_2 \)
  \eq
  [\alpha]_q^{1/2} \ket{3} \bra{1}
  -
  {[\alpha +1 ]}_q^{1/2} \ket{4} \bra{2}.
\ee

\pagebreak

Associated with $U_q\[gl(2|1)\]$ there is  a co-product structure
($\mathbb{Z}_2$-graded algebra homomorphism)
$
  \Delta : U_q\[gl(2|1)\] \to U_q\[gl(2|1)\] \otimes U_q\[gl(2|1)\]
$
given by:
\be
  \Delta ( {E^i}_i)
  \eq
  I\otimes {E^i}_i + {E^i}_i \otimes I,
  \qquad
  i = 1, 2, 3,
  \\
  \Delta ( {E^1}_2 )
  \eq
  {E^1}_2 \otimes q^{-\frac{1}{2} ( {E^1}_1 - {E^2}_2 )}
  +
  q^{\frac{1}{2} ( {E^1}_1 - {E^2}_2 )} \otimes {E^1}_2
  \\
  \Delta ( {E^2}_1 )
  \eq
  {E^2}_1 \otimes q^{-\frac{1}{2} ( {E^1}_1 - {E^2}_2)}
  +
   q^{\frac{1}{2} ( {E^1}_1 - {E^2}_2 )} \otimes {E^2}_1
  \\
  \Delta ( {E^2}_3 )
  \eq
  {E^2}_3 \otimes q^{-\frac{1}{2} ( {E^2}_2 + {E^3}_3 )}
  +
  q^{\frac{1}{2}( {E^2}_2 + {E^3}_3 )} \otimes {E^2}_3
  \\
  \Delta ( {E^3}_2 )
  \eq
  {E^3}_2 \otimes q^{-\frac{1}{2}( {E^2}_2 + {E^3}_3 )}
  +
  q^{\frac{1}{2}( {E^2}_2 + {E^3}_3 )} \otimes {E^3}_2.
\ee
There exists another possible co-product structure:
$ \overline{\Delta} $, defined by
$ \overline{\Delta} = T \cdot \Delta $, where
$
  T :
  U_q \[ gl \( 2 | 1 \) \] \otimes U_q \[ gl \( 2 | 1 \) \]
  \to
  U_q \[ gl \( 2 | 1 \) \] \otimes U_q \[ gl \( 2 | 1 \) \]
$
is the twist map, defined for homogeneous elements
$ a, b \in U_q \[ gl \( 2 | 1 \) \] $:
\be
  T \( a \otimes b \)
  =
  {\( - \)}^{\[ a \] \[ b \]}
  \( b \otimes a \).
\ee

The tensor product module has the following decomposition with respect
to the co-product for generic values of $ \alpha $:
\bse
  V \otimes V
  =
  V_1 \oplus V_2 \oplus V_3.
  \eqlabel{TPdecomp}
\ese
We construct symmetry adapted bases
$ \{ \ket{\Psi^k_1} \}_{k=1}^4 $ and $ \{ \ket{\Psi^k_3} \}_{k=1}^4 $,
for the spaces $ V_1 $ and $ V_3 $ respectively in terms of the basis
elements of $ V $:
\be
  \ket{\Psi^1_1}
  \eq
  \ket{1} \otimes \ket{1}
  \\
  \ket{\Psi^1_2}
  \eq
  {(q^{\alpha}+q^{-\alpha})}^{-\frac{1}{2}}
  \(
    q^{\alpha/2} \ket{1} \otimes \ket{2}
    +
    q^{-\alpha/2} \ket{2} \otimes \ket{1}
  \)
  \\
  \ket{\Psi^1_3}
  \eq
  (q^{\alpha}+q^{-\alpha})^{-\frac{1}{2}}
  \(
    q^{\alpha/2} \ket{1} \otimes \ket{3}
    +
    q^{-\alpha/2} \ket{3} \otimes \ket{1}
  \)
  \\
  \ket{\Psi^1_4}
  \eq
  {\( q^{\alpha} + q^{-\alpha} \)}^{-\frac{1}{2}}
  {[2\alpha+1]}_q^{-\frac{1}{2}}
  \times
  \\
  & &
  \! \! \! \! \! \! \! \!
  \! \! \! \! \! \! \! \!
  \[
    {[\alpha+1]}_q^{\frac{1}{2}}
    \(
      q^{\alpha} \ket{1} \otimes \ket{4}
      +
      q^{-\alpha} \ket{4} \otimes \ket{1}
    \)
    +
    {[\alpha]}_q^{\frac{1}{2}}
    \(
      q^{\frac{1}{2}} \ket{2} \otimes \ket{3}
      -
      q^{-\frac{1}{2}} \ket{3} \otimes \ket{2}
    \)
  \]
  \\
  \ket{\Psi^3_1}
  \eq
  {(q^{\alpha+1}+q^{-\alpha-1})}^{-\frac{1}{2}}
  {[2\alpha+1]}_q^{-\frac{1}{2}}
  \times
  \\
  & &
  \! \! \! \!
  \! \! \! \!
  \! \! \! \!
  \! \! \! \!
  \! \! \! \!
  \! \! \! \!
  \[
    {[\alpha]}_q^{\frac{1}{2}}
    \(
      q^{\alpha+1} \ket{4} \otimes \ket{1}
      +
      q^{-\alpha-1} \ket{1} \otimes \ket{4}
    \)
    +
    [\alpha+1]_q^{\frac{1}{2}}
    \(
      q^{-\frac{1}{2}} \ket{3} \otimes \ket{2}
      -
      q^{\frac{1}{2}} \ket{2} \otimes \ket{3}
    \)
  \]
  \\
  \ket{\Psi^3_2}
  \eq
  {(q^{\alpha+1}+q^{-\alpha-1})}^{-\frac{1}{2}}
  \(
    q^{(\alpha+1)/2} \ket{4} \otimes \ket{2}
    +
    q^{-(\alpha+1)/2} \ket{2} \otimes \ket{4}
  \)
  \\
  \ket{\Psi^3_3}
  \eq
  {(q^{\alpha+1}+q^{-\alpha-1})}^{-\frac{1}{2}}
  \(
    q^{(\alpha+1)/2} \ket{4} \otimes \ket{3}
    +
    q^{-(\alpha+1)/2} \ket{3} \otimes \ket{4}
  \)
  \\
   \ket{\Psi^3_4}
  \eq
   \ket{4} \otimes \ket{4}.
\ee

Dual bases
$ \{ \bra{\Psi^k_1} \}_{k=1}^4 $ and $ \{ \bra{\Psi^k_3} \}_{k=1}^4 $,
are found from the definitions:
\bne
  \bra{\Psi^k_j}
  \eq
  \ket{\Psi^k_j}^{\dagger},
  \qquad \qquad \qquad
  k = 1, 3,
  \qquad
  j = 1, \dots, 4,
  \eqlabel{dualrule1}
  \\
  {( \ket{i} \otimes \ket{j} )}^{\dagger}
  \eq
  {\( - \)}^{\[ \ket{i} \] \[ \ket{j} \]}
  \( \bra{i} \otimes \bra{j} \),
  \qquad
  i, j = 1, \dots, 4.
  \eqlabel{dualrule2}
\ene
Now, the general form of the basis vectors $ \ket{\Psi^k_j} $ is:
\be
  \ket{\Psi^k_j}
  =
  \sum_{m}
    \theta^{k j}_{m}
    \( \ket{x^{k j}_{m}} \otimes \ket{y^{k j}_{m}} \),
\ee
where the $ \theta^{k j}_{m} $ are in general complex scalar
functions of $ q $ and $ \alpha $. From \eqref{dualrule1} and
\eqref{dualrule2}, and choosing the parameters $ q $ and $ \alpha $ to
be real and positive, the duals of these vectors are given by:
\be
  \bra{\Psi^k_j}
  =
  \sum_{m}
    {\( - \)}^{ \[ \ket{x^{k j}_{m}} \] \[ \ket{y^{k j}_{m}} \]}
    \theta^{k j}_{m}
    \( \bra{x^{k j}_{m}} \otimes \bra{y^{k j}_{m}} \).
\ee
As the $ R $ matrix is unique,
analytic continuation makes our final results valid for any
complex $ q $ and $ \alpha $. For the duals, we obtain:
\be
  \bra{\Psi^1_1}
  \eq
  \bra{1} \otimes \bra{1}
  \\
  \bra{\Psi^1_2}
  \eq
  (q^{\alpha}+q^{-\alpha})^{-\frac{1}{2}}
  \(
    q^{\frac{1}{2}\alpha} \bra{1} \otimes \bra{2}
    +
    q^{-\frac{1}{2}\alpha} \bra{2} \otimes \bra{1}
  \)
  \\
  \bra{\Psi^1_3}
  \eq
  (q^{\alpha}+q^{-\alpha})^{-\frac{1}{2}}
  \(
    q^{\frac{1}{2}\alpha} \bra{1} \otimes \bra{3}
    +
    q^{-\frac{1}{2}\alpha} \bra{3} \otimes \bra{1}
  \)
  \\
  \bra{\Psi^1_4}
  \eq
  {\( q^{\alpha} + q^{-\alpha} \)}^{-\frac{1}{2}}
  [2\alpha+1]_q^{-\frac{1}{2}}
  \times
  \\
  & &
  \! \! \! \!
  \! \! \! \!
  \! \! \! \!
  \[
    [\alpha+1]_q^{\frac{1}{2}}
    \(
      q^{\alpha} \bra{1} \otimes \bra{4}
      +
      q^{-\alpha} \bra{4} \otimes \bra{1}
    \)
    -
    [\alpha]_q^{\frac{1}{2}}
    \(
      q^{\frac{1}{2}} \bra{2} \otimes \bra{3}
      -
      q^{-\frac{1}{2}} \bra{3} \otimes \bra{2}
    \)
  \]
  \\
  \bra{\Psi^3_1}
  \eq
  (q^{\alpha+1}+q^{-\alpha-1})^{-\frac{1}{2}}
  [2\alpha+1]_q^{-\frac{1}{2}}
  \times
  \\
  & &
  \! \! \! \!
  \! \! \! \!
  \! \! \! \!
  \! \! \! \!
  \! \! \! \!
  \! \! \! \!
  \[
    [\alpha]_q^{\frac{1}{2}}
    \(
      q^{\alpha+1} \bra{4} \otimes \bra{1}
      +
      q^{-\alpha-1} \bra{1} \otimes \bra{4}
    \)
    -
    [\alpha+1]_q^{\frac{1}{2}}
    \(
      q^{-\frac{1}{2}} \bra{3} \otimes \bra{2}
      -
      q^{\frac{1}{2}} \bra{2} \otimes \bra{3}
    \)
  \]
  \\
  \bra{\Psi^3_2}
  \eq
  {\( q^{\alpha+1}+q^{-\alpha-1} \)}^{-\frac{1}{2}}
  \(
    q^{\frac{1}{2}(\alpha+1)} \bra{4} \otimes \bra{2}
    +
    q^{-\frac{1}{2}(\alpha+1)} \bra{2} \otimes \bra{4}
  \)
  \\
  \bra{\Psi^3_3}
  \eq
  (q^{\alpha+1}+q^{-\alpha-1})^{-\frac{1}{2}}
  \(
    q^{\frac{1}{2}(\alpha+1)} \bra{4} \otimes \bra{3}
    +
    q^{-\frac{1}{2}(\alpha+1)} \bra{3} \otimes \bra{4}
  \)
  \\
  \bra{\Psi^3_4}
  \eq
  \bra{4} \otimes \bra{4}.
\ee

From the basis vectors $ \ket{\Psi^k_j} $ and their duals
$ \bra{\Psi^k_j} $ for $ V_1 $ and $ V_3 $, we
construct projectors $ P_1 $ and $ P_3 $, defined by:
\be
  P_1
  =
  \sum_{k=1}^4
    \ket{\Psi^k_1} \bra{\Psi^k_1},
  \qquad
  P_3
  =
  \sum_{k=1}^4
    \ket{\Psi^k_3} \bra{\Psi^k_3}.
\ee
Note that the multiplication operation on the graded space
$ V \otimes V $ is given by:
\bse
  \( \ket{i} \otimes \ket{j} \) \( \bra{k} \otimes \bra{l} \)
  =
  {\( - \)}^{\[ \ket{j} \] \[ \bra{k} \]}
  \( \ket{i} \bra{k} \otimes \ket{j} \bra{l} \),
  \qquad
  i, j, k, l = 1, 2, 3, 4.
  \eqlabel{Callistemon}
\ese
Now let $ I $ be the identity operator on $ V \otimes V $, viz:
$
  I
  =
  \sum_{i j=1}^4
    {e^i}_i \otimes {e^j}_j
$,
where $ {e^k}_l = \ket{k} \bra{l} $ is an elementary rank $ 2 $ tensor.
As we have $ P_1 + P_2 + P_3 = I $, we thus do not need to explicitly
construct $ P_2 $ (or even a basis for $ V_2 $); we simply set:
\bse
  P_2 = I - P_1 - P_3.
  \eqlabel{P2def}
\ese


Where $g$ is a classical Lie superalgebra, the corresponding quantum
superalgebra $U_q[g]$ admits a universal $R$ matrix
$R\in U_q[g]\otimes U_q[g]$ satisfying (among other relations):
\bne
  R \Delta \( a \)
  \eq
  \overline{\Delta} \( a \) R,
  \qquad \quad \;
  \forall a \in U_q \[ g \],
  \nonumber
  \\
  R_{12} R_{13} R_{23}
  \eq
  R_{23} R_{13} R_{12},
  \qquad
  \mathrm{in}
  \;
  U_q \[ g \] \otimes U_q \[ g \] \otimes U_q \[ g \],
  \eqlabel{YB}
\ene
where the subscripts refer to the embedding of $ R $ acting on the
triple tensor product space. From any representation of
$ U_q \[ g \] $, one may obtain a tensor solution of \eqref{YB} by
replacing the superalgebra elements with their matrix representatives.
Similarly to \eqref{Callistemon}, multiplication of tensor products of
matrices $ a, b, c, d $ is governed by:
\be
  \( a \otimes b \)
  \( c \otimes d \)
  =
  {\( - \)}^{\[ b \] \[ c \]}
  \( a c \otimes b d \),
  \qquad
  \textrm{homogeneous} \;
  b, c.
\ee

We introduce the \emph{graded permutation operator $ P $} on the
tensor product space $ V \otimes V $,
defined for graded basis vectors
$ v^k, v^l \in V $ by:
\be
  P ( v^k \otimes v^l )
  =
  {\( - \)}^{\[ k \] \[ l \]}
  ( v^l \otimes v^k ),
\ee
and extended by linearity.
(We use the shorthand $ \[ v^k \] \equiv \[ k \] $.)
With this, we define:
\be
  \sigma = P R,
\ee
which can be shown to satisfy the equation:
\bse
  \( \sigma \otimes I \)
  \( I \otimes \sigma \)
  \( \sigma \otimes I \)
  =
  \( I \otimes \sigma \)
  \( \sigma \otimes I \)
  \( I \otimes \sigma \).
  \eqlabel{Isigma}
\ese

From \cite{LinksGould:92b}, we have (with a slight change of notation
and a convenient choice of normalisation):
\be
  \sigma
  =
  q^{-2 \alpha} P_1
  -
  P_2
  +
  q^{2 \alpha + 2} P_3.
\ee
Using \eqref{P2def}, this simplifies to:
\be
  \sigma
  =
  \( 1 + q^{-2 \alpha} \) P_1
  +
  \( 1 + q^{2 \alpha + 2} \) P_3
  -
  I.
\ee

From the above form of $ \sigma $, it is straightforward to deduce that
$ \sigma $ satisfies the polynomial identity:
\be
  q^{-1} \sigma^3
  +
  \( q^{-1} - q^{- 2 \alpha - 1} - q^{2 \alpha + 1} \) \sigma^2
  +
  \( q - q^{- 2 \alpha - 1} - q^{2 \alpha + 1} \) \sigma
  +
  q I
  =
  0.
\ee
The above skein relation may be used to evaluate the invariant in some
cases, but not all since it is of third order.  The invariant may also
be directly evaluated for a class of links using quantum superalgebra
theoretic results \cite{GouldLinksZhang:96b}.

\pagebreak

We will represent rank $ 2 $ tensors as matrices, that is, the
elementary rank $ 2 $ tensor $ {e^i}_k $ is represented by the
elementary ($ 4 \times 4 $) matrix $ e_{i,k} $.
We adopt the (standard) convention that the elementary rank $ 4 $
tensor $ {e^{i j}}_{k l} = {e^i}_k \otimes {e^j}_l $ is constructed by
insertion of a copy of the elementary rank $ 2 $ tensor $ {e^j}_l $ at
each location of $ {e^i}_k $ (i.e.
each element of $ {e^i}_k $ is multiplied by the whole of $ {e^j}_l $).
This means that
$ {e^{i j}}_{k l} $ is represented by the elementary ($ 16 \times 16 $)
matrix $ e_{4(i-1)+j,4(k-1)+l} $.

Let $ A $ be an arbitrary graded rank $ 4 $ tensor acting on
$ V \otimes V $, then for scalar coefficients $ {A^{i j}}_{k l} $:
\be
  A
  =
  \sum_{i j k l}
    {A^{i j}}_{k l}
    \( {e^i}_k \otimes {e^j}_l \).
\ee
Our convention then tells us that the coefficient $ {A^{i j}}_{k l} $ is
the $ \( 4(i-1)+j, 4(k-1)+l \) $ entry of
$ A $, written explicitly:
\be
  {A^{i j}}_{k l}
  \mapsto
  A_{4(i-1)+j,4(k-1)+l}.
\ee

We wish to remove the grading on $ V $, and convert the matrix
representing $ \sigma $ to its ungraded counterpart.
Recall that basis vectors $ v^k $ satisfy
$ {e^i}_j v^k = {\delta^k}_j v^i $, hence the action of $ A $ on
$ V \otimes V $ is:
\be
  A ( v^k \otimes v^l )
  \eq
  \sum_{i j m n}
    {A^{i j}}_{m n}
    \( {e^i}_m \otimes {e^j}_n \)
    ( v^k \otimes v^l )
  \\
  \eq
  \sum_{i j m n}
    {A^{i j}}_{m n}
    {\( - \)}^{\[ k \] \( \[ j \] + \[ n \] \)}
    \( {e^i}_m v^k \otimes {e^j}_n v^l \)
  \\
  \eq
  \sum_{i j m n}
    {A^{i j}}_{m n}
    {\( - \)}^{\[ k \] \( \[ j \] + \[ n \] \)}
    \( {\delta^k}_m v^i \otimes {\delta^l}_n v^j \)
  \\
  \eq
  \sum_{i j m n}
    {A^{i j}}_{m n}
    {\( - \)}^{\[ k \] \( \[ j \] + \[ n \] \)}
    {\delta^k}_m {\delta^l}_n
    \( v^i \otimes v^j \)
  \\
  \eq
  \sum_{i j}
    {A^{i j}}_{k l}
    {\( - \)}^{\[ k \] \( \[ j \] + \[ l \] \)}
    \( v^i \otimes v^j \).
\ee
Now, in this sum, the parity factor is constructed from the degrees
of vectors; in the ungraded case, there would be no such factor, indeed
we would have:
\be
  \overline{A} ( v^k \otimes v^l )
  \eq
  \sum_{i j}
    {\overline{A}^{i j}}_{k l}
    \( v^i \otimes v^j \).
\ee
This motivates us to set:
\be
  {\overline{A}^{i j}}_{k l}
  =
  {\( - \)}^{\[ k \] \( \[ j \] + \[ l \] \)}
  {A^{i j}}_{k l}.
\ee

Under these conventions, the explicit form of $ \sigma $ is presented
(as a matrix!) in \secref{lglink}.


\section{Knot Theory}

%
%

\subsection{Link Examples}

In \tabref{linknamesandwrithes}, we list the links to be studied.  (By
the term `knot', we intend a link of one component.) We use the
well-known notation of Alexander and Briggs (1926)
\cite{AlexanderBriggs:26}, the data being abstracted from
\cite{Adams:94}, itself citing \cite{Rolfsen:76} and
\cite{DollHoste:91} (beware that the tables in this latter article are
presented in microfiche form only).

\def\cross{No}
\def\tick{Yes}
\def\scs{\scriptstyle}

\begin{table}[ht]
  \centering
  \begin{tabular}{||l|c|p{32mm}|p{32mm}||}
    \hline\hline
    & & & \\[-3mm]
      \multicolumn{1}{||c|}{$ K $} &
      $ w \( K \) $ &
      \multicolumn{1}{c|}{Chiral?} &
      \multicolumn{1}{c||}{Invertible?} \\
    & & & \\[-3mm]
    \hline\hline
    & & & \\[-2mm]
    $ 0_1 $ (Unknot) &
      $ 0 $ & $ \cross $ (trivial) & $ \tick $ (trivial) \\[2mm]
    \hline
    & & & \\[-2mm]
    $ 2^2_1 $ (Hopf Link) &
      $ 2 $ & $ \cross $ (trivial) & $ \tick $ (trivial) \\[2mm]
    \hline
    & & & \\[-2mm]
      
      $ 3_1 $ (Trefoil) &
      $ 3 $ & $ \tick $ \cite[p~176]{Adams:94} & $ \tick $ (trivial)
      \\[2mm]
    \hline
    & & & \\[-2mm]
      $ 4_1 $ (Figure Eight) &
      $ 0 $ &
      $ \cross $
        (\cite[p~14]{Adams:94};
        see \cite[p~198]{Kauffman:88} for
        an elegant graphical proof) &
      $ \tick $ ({\raggedright as $ 8_{17} $ is \\}
      the smallest noninvertible knot)
      \\[2mm]
    \hline
    & & & \\[-2mm]
    $ 5^2_1 $ (Whitehead Link) &
      $ 1 $ & $ \tick $ \cite[pp~49-50]{Kauffman:87a} & $ \tick $
      \\[2mm]
    \hline
    & & & \\[-2mm]
    $ 8_{17} $ &
      $ 0 $ & $ \cross $ \cite[p~455]{Kauffman:87a} &
      $ \cross $ \cite[p~162]{Kauffman:97a}
      \\[2mm]
    \hline
    & & & \\[-2mm]
    $ 9_{42} $ &
      $ 1 $ & $ \tick $ \cite[p~218]{Kauffman:93} & $ \tick $ \\[2mm]
    \hline
    & & & \\[-2mm]
    $ 10_{48} $ &
      $ 0 $ & $ \tick $ \cite[p~218]{Kauffman:93} & $ \tick $ \\[2mm]
    \hline\hline
  \end{tabular}
  \caption{%
    Data for the links to be investigated, including their
    Alexander--Briggs (and common) names,
    their writhes $ w \( K \) $, and whether they are chiral and
    invertible.  Diagrams of the links are presented in Figures
    \ref{fig:HopfLinkTrefoil} to \ref{fig:TenFortyEight}.
  }
  \tablabel{linknamesandwrithes}
\end{table}


\subsection{Reflection and Inversion -- Chirality and Invertibility}

Throughout, we shall write ``='' to denote \emph{ambient isotopy} of
link diagrams, meaning that they are equivalent under the Reidemeister
moves (original: \cite{Reidemeister:48}, but see, e.g.
\cite{Kauffman:93}).  We shall use the following definitions, but the
reader must be aware that conflicting terminology appears in the
literature.

\pagebreak

\begin{description}
\item[Reflection:]
  We shall denote by $ K^* $ the mirror image (or reflection) of a knot
  $ K $.  A knot is \emph{chiral} if it is distinct from its mirror
  image; i.e. there are actually two distinct knots with the same name,
  $ K^* \neq K $, e.g. the trefoil knot is chiral:  $ {\( 3_1 \)}^*
  \neq 3_1 $.  Note that this definition doesn't require an
  orientation. A knot is \emph{amphichiral} if it is ambient isotopic
  to its mirror image, i.e. $ K^* = K $.

  The HOMFLY%
  \footnote{
    The HOMFLY polynomial is named by the conjunction of the initials
    of six of its discoverers
    \cite{FreydYetterHosteLickorishMilletOcneanu:85}, omitting those
    (``P'' and ``T'') of
    two independent discoverers \cite{PrzytyckiTraczyk:87}.
    Przytycki, the omitted ``P'', has furthered the entymological
    spirit with the suggestion ``FLYPMOTH'' \cite[p~256]{Przytycki:89},
    which includes all the discoverers and has a muted reference to the
    ``flyping'' operation of the Tait, Kirkwood and Little -- the
    original compilers of knot tables.
    (Another possibility is the letter sequence ``HOMFLYPT''.)
    Bar-Natan (Prasolov and Sossinsky \cite[p~36]{PrasolovSossinsky:96}
    cite Bar-Natan \cite{BarNatan:95}, who cites ``L Rudulph'')
    goes further, adding a ``U'' for good measure, to account for any
    unknown discoverers, yielding the unpalatable ``LYMPHTOFU''!
  }
 (and hence the Jones) polynomial and the
  Kauffman polynomial can distinguish many (but not all) knots from
  their reflections.
  The first chiral knot that neither the HOMFLY nor the Kauffman
  polynomial can distinguish is $ 9_{42} $,
  i.e. $ 9_{42}^* \neq 9_{42} $, but the polynomials are equal.
  Similarly, the knot $ 10_{48} $ is chiral, but the HOMFLY polynomial
  fails to detect this, although the Kauffman does detect it
  \cite[p~218]{Kauffman:93} (wrongly labeled $10_{79}$).
\item[Inversion:]
  Assign an orientation to a knot. Denote the \emph{inverse}
  of a knot $ K $ by $ K^{-1} $, obtained by reversing
  the orientation.  Whilst this is a simple concept for a knot, there
  are of course many possibilities for the reversal of only some
  components of oriented, multi-component links; we shall not go into
  these here.

  Commonly, $ K = K^{-1} $, and we say that $ K $ is
  \emph{invertible}.  For example, the trefoil knot is invertible
  $ {\( 3_1 \)}^{-1} = 3_1 $.  Less commonly, $ K \neq K^{-1} $, and we
  say that $ K $ is \emph{noninvertible}.  The first example of a
  noninvertible (prime) knot is $ 8_{17} $.
\end{description}

Both the reflection and the inverse are automorphisms of order two, i.e.
$ {( K^* )}^* = K $ and $ {( K^{-1} )}^{-1} = K $.
The notions may of course be combined, we obtain:
$ {( K^* )}^{-1} = {( K^{-1} )}^* $.

To illustrate, using the Trefoil Knot $ 3_1 $
(see \figref{HopfLinkTrefoil}). We have two equivalence classes:
$ 3_1 = {\( 3_1 \)}^{-1} $ and
$ {\( 3_1 \)}^* = {( {\( 3_1 \)}^{-1} )}^* = {( \( 3_1 \) ^* )}^{-1} $.


\subsection{Abstract Tensor Conventions}

By a `positive oriented' or `right-handed' crossing, we shall intend a
crossing such that if the thumb of the right hand points in the
direction of one of the arrows, the fingers of the right hand will
point in the direction of the other arrow. The opposite situation is
naturally called a `negative oriented' or `left-handed' crossing.

If the two outward-pointing arrows of a positive oriented crossing are
pointed upwards, then we shall label the components of the crossing
with indices $ a $ in top left, $ b $ in bottom left, $ c $ in top
right, and $ d $ in bottom right, and associate with the crossing the
(rank $ 4 $) tensor $ {\sigma}^{a~c}_{b~d} $, where the position of the
indices in the tensor corresponds with the positioning of the labels in
the crossing.  The inverse of $ \sigma $ will represent a negative
oriented crossing, with the convention on the indices being the same as
that of $ \sigma $.  A diagram of $ \sigma $ and $ \sigma^{-1} $ is
provided in \figref{RandS}.

\begin{figure}[htbp]
  \begin{center}
    \input{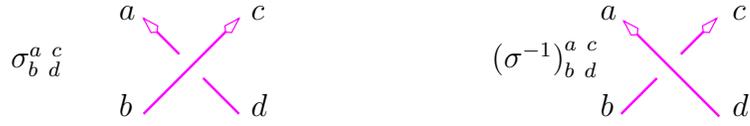}
    \caption{
      Definition of the tensors $ \sigma $ and $ \sigma^{-1} $
      representing positive oriented and negative oriented crossings
      with upward pointing arrows, respectively.
    }
    \figlabel{RandS}
  \end{center}
\end{figure}

We shall also require four (rank $ 2 $) tensors (i.e. genuine matrices)
to represent all possible horizontally-oriented half-loops. We shall
call these `cap' and `cup' matrices, and label them with the suggestive
$ \Omega^\pm $ and $ \mho^\pm $, e.g. $ \Omega^+ $ is the upper loop
with arrow pointing right.  A diagram is provided in
\figref{CapsandCups}.

\begin{figure}[htbp]
  \begin{center}
    \input{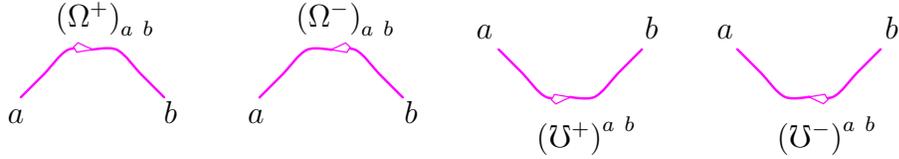}
    \caption{
      Definition of the tensors (matrices)
      $ \Omega^\pm $ and $ \mho^\pm $, representing all possibilities
      of horizontally-aligned half-loops.
    }
    \figlabel{CapsandCups}
  \end{center}
\end{figure}

With these basic tensors $ \sigma $, $ \sigma^{-1} $, $ \Omega^\pm $
and $ \mho^\pm $, we may evaluate an invariant for any particular
link.  However, this direct procedure tends to be computationally
expensive, and parts of the computation are often repeated, so in
practice, we define auxiliary symbols.  We shall use the notation
$ X $ to represent a rank $ 4 $ tensor such as $ \sigma $ or
$ \sigma^{-1} $ with parallel pointing arrows
(i.e. a `channel' crossing in the terminology of Kauffman
\cite[p~76]{Kauffman:93}.)

The primary auxiliary tensors used are listed below; secondary ones
will be mentioned where necessary.  The Einstein summation convention
is used throughout.

\begin{itemize}
\item
  The first auxiliary symbols are those of crossings that have been
  `twisted' relative to $ \sigma $ and $ \sigma^{-1} $. The left,
  right, and upside-down-twisted versions of $ X $ will be called
  $ X_l $, $ X_r $ and $ X_d $ respectively. They are defined in the
  following manner:
  \bne
    {\( X_l \)}^{a~c}_{b~d}
    & \defeq &
    {X}^{e~a}_{d~h}
    \cdot
    {\( \Omega^- \)}_{b~e}
    \cdot
    {\( \mho^- \)}^{h~c}
    \nonumber
    \\
    {\( X_r \)}^{a~c}_{b~d}
    & \defeq &
    {X}^{c~g}_{f~b}
    \cdot
    {\( \mho^+ \)}^{a~f}
    \cdot
    {\( \Omega^+ \)}_{g~d}
    \eqlabel{XRdef}
    \\
    {\( X_d \)}^{a~c}_{b~d}
    & \defeq &
    {X}^{e~g}_{f~h}
    \cdot
    {\( \mho^+ \)}^{a~h}
    \cdot
    {\( \Omega^+ \)}_{g~b}
    \cdot
    {\( \mho^+ \)}^{c~f}
    \cdot
    {\( \Omega^+ \)}_{e~d}.
    \nonumber
  \ene
  Observe that $ X_d $ is a `channel' crossing, whilst $ X_l $ and
  $X_r$ are `cross-channel' crossings.  Diagrams are found in Figures
  \ref{fig:XLandXR} and \ref{fig:XD}.

  \begin{figure}[htbp]
    \begin{center}
      \input{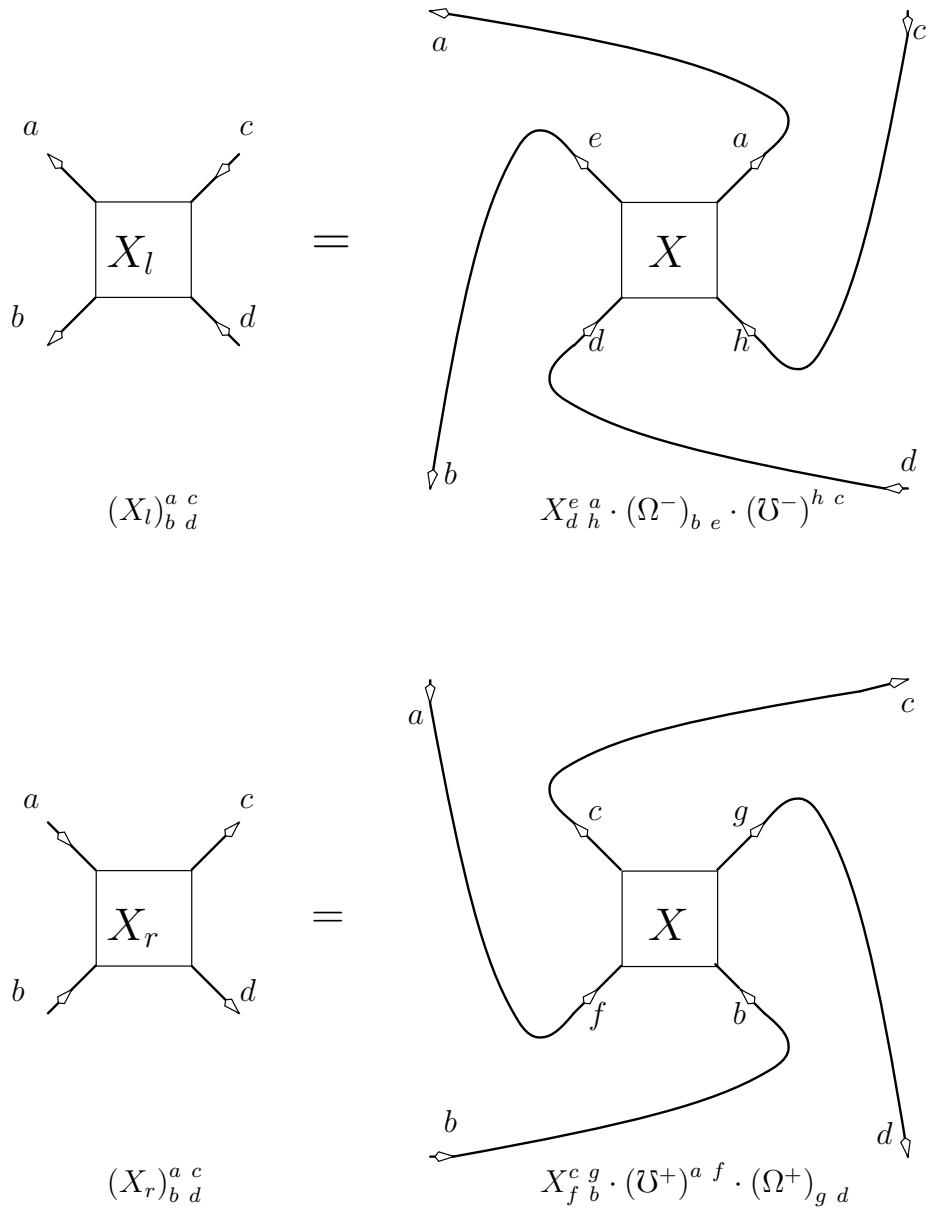}
      \caption{The primary auxiliary tensors
        $ X_l $ and $ X_r $, where
        $ X $ is one of $ \sigma $ or $ \sigma^{-1} $.
      }
      \figlabel{XLandXR}
    \end{center}
  \end{figure}

  \begin{figure}[htbp]
    \begin{center}
      \input{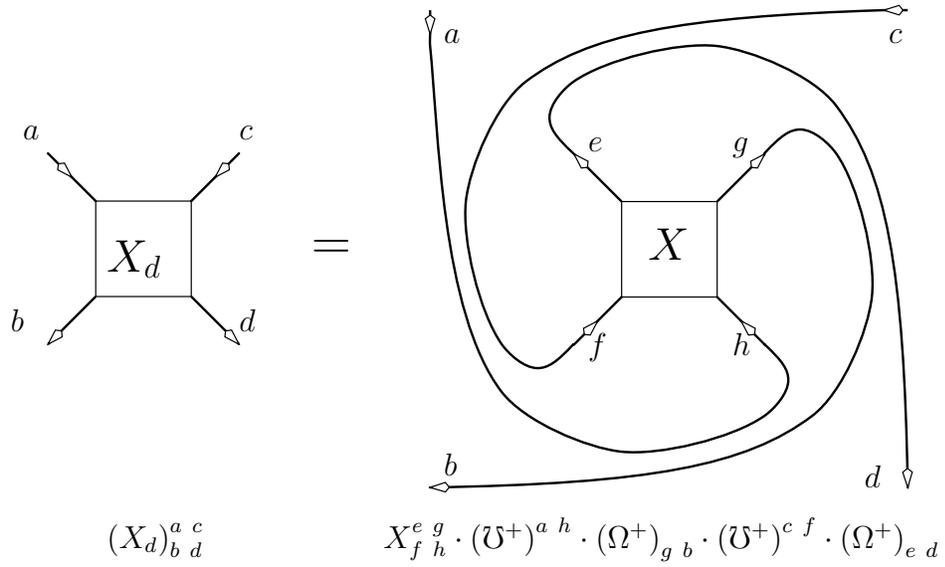}
      \caption{The primary auxiliary tensor
        $ X_d $, where
        $ X $ is one of $ \sigma $ or $ \sigma^{-1} $.
      }
      \figlabel{XD}
    \end{center}
  \end{figure}
\item
  The next set of auxiliary symbols represent $ p $ copies of the
  \emph{same} crossing $ X $ (for any channel crossing $ X $) atop one
  another (see \figref{Xp}).  They are defined recursively in the
  following manner:
  \be
    {\( X^{p+1} \)}^{a~c}_{b~d}
    \defeq
    X^{a~c}_{e~f}
    \cdot
    {\( X^p \)}^{e~f}_{b~d},
    \qquad
    p = 1, 2, \dots.
  \ee

  \begin{figure}[htbp]
    \begin{center}
      \input{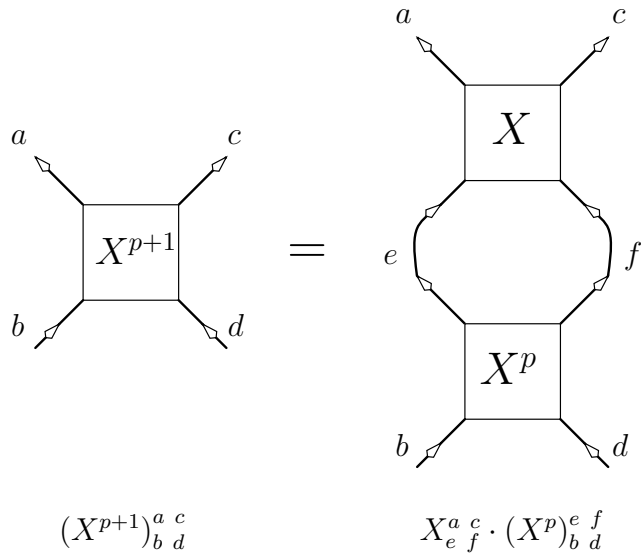}
      \caption{The primary auxiliary tensors
        $ X^{p+1} $ in terms of $ X $ and $ X^p $;
	$ X $ is one of $ \sigma $ or $ \sigma^{-1} $. If all
	arrows are reversed, then the definition also holds for $ X $
	being $ \sigma_d $ or $ \sigma^{-1}_d $; that is, any channel
	crossing.
      }
      \figlabel{Xp}
    \end{center}
  \end{figure}
\item
  The third set of frequently-encountered patterns are where a crossing
  $ X $ is to the left or right of its own `upside-downness' $ X_d $.
  That is, fix an $ X $ as either $ \sigma $ or $ \sigma^{-1} $, and
  examine the patterns formed from juxtaposing $ X $ and $ X_d $.  They
  are defined in the following manner:
  \be
    {\( X_d X \)}^{a~c}_{b~d}
    & \defeq &
    {\( X_d \)}^{a~e}_{b~f}
    \cdot
    {X}^{g~c}_{h~d}
    \cdot
    {\( \Omega^- \)}_{e~g}
    \cdot
    {\( \mho^+ \)}^{f~h}
    \\
    {\( X X_d \)}^{a~c}_{b~d}
    & \defeq &
    {X}^{a~e}_{b~f}
    \cdot
    {\( X_d \)}^{g~c}_{h~d}
    \cdot
    {\( \Omega^+ \)}_{e~g}
    \cdot
    {\( \mho^- \)}^{f~h}.
  \ee
  A diagram is found in \figref{XXDandXDX}.

  \begin{figure}[htbp]
    \begin{center}
      \input{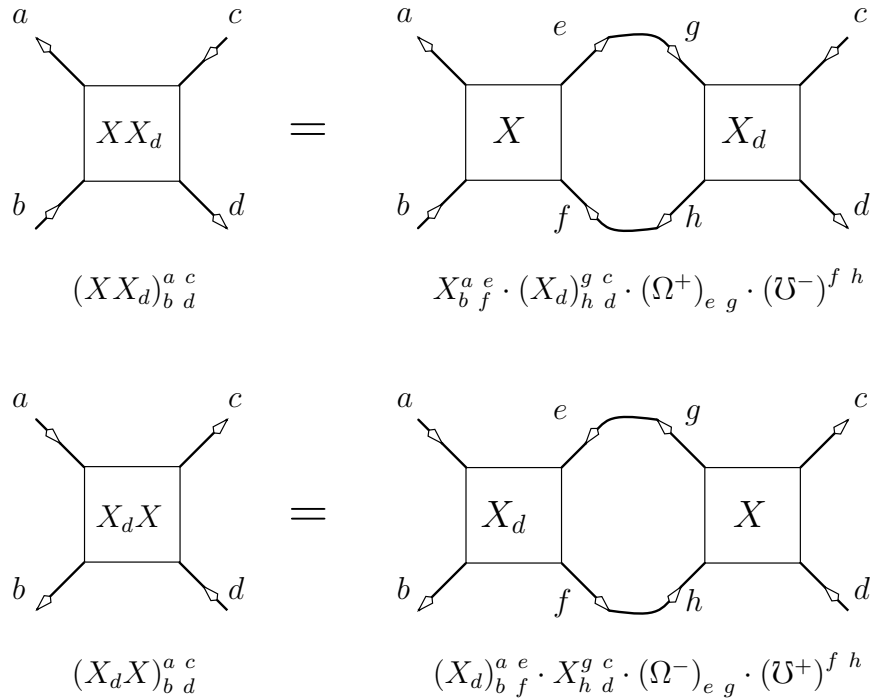}
      \caption{The primary auxiliary tensors
        $ X_d X $ and $ X X_d $;
        $ X $ is one of $ \sigma $ or $ \sigma^{-1} $.
      }
      \figlabel{XXDandXDX}
    \end{center}
  \end{figure}
\item
  The final set of frequently-encountered patterns are where a crossing
  $ X_l $ is placed atop above a crossing $ X_r $ (or vice-versa). We
  obtain:
  \be
    {\( X_l X_r \)}^{a~c}_{b~d}
    & \defeq &
    {\( X_l \)}^{a~c}_{e~f}
    \cdot
    {\( X_r \)}^{e~f}_{b~d}
    \\
    {\( X_r X_l \)}^{a~c}_{b~d}
    & \defeq &
    {\( X_r \)}^{a~c}_{e~f}
    \cdot
    {\( X_l \)}^{e~f}_{b~d}.
  \ee
  A diagram is found in \figref{XLXRandXRXL}.  A moment's thought
  demonstrates that the diagram for $ X_l X_r $ is a right rotation of
  the diagram for $ X_d X $. In fact, we have the identity:
  \be
    {\( X_l X_r \)}^{a~c}_{b~d}
    =
    {\( X_d X \)}^{e~a}_{d~h}
    \cdot
    {\( \mho^+ \)}^{h~c}
    \cdot
    {\( \Omega^+ \)}_{b~e},
  \ee
  although in practice we shall not use it. (A diagram parallel to
  \figref{XLandXR} would demonstrate this.)

  \begin{figure}[htbp]
    \begin{center}
      \input{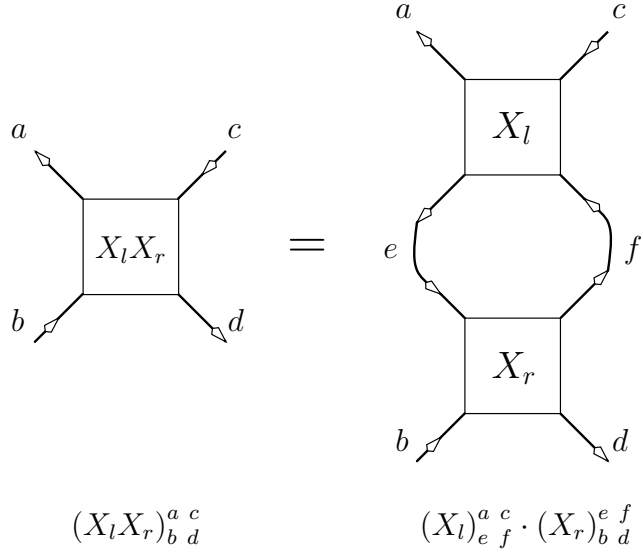}
      \caption{The primary auxiliary tensor
        $ X_l X_r $;
        $ X $ is one of $ \sigma $ or $ \sigma^{-1} $.
        $ X_r X_l $ is obtained by swapping every $ r $ and $ l $
        in this diagram.
      }
      \figlabel{XLXRandXRXL}
    \end{center}
  \end{figure}
\end{itemize}


\subsection{The Effects of Reflection and Inversion on the Tensors}

\begin{description}
\item[Reflection:]
  Let $ K^* $ be the reflection of a tangle $ K $; and say that we have
  constructed a tensor representing $ K $. Every positive (respectively
  negative) crossing in $ K $ will have been replaced by the equivalent
  negative (respectively positive) crossing in $ K^* $.  Thus, the
  tensor corresponding to $ K^* $ will be that of $ K $ with every
  $ \sigma $ replaced by $ \sigma^{-1} $, and every $ \sigma^{-1} $
  replaced by $ \sigma $. This carries through to the auxiliary
  tensors; i.e.  $ \sigma_d \sigma $ will be replaced with
  $ \sigma^{-1}_d \sigma^{-1} $, etc.  The caps $ \Omega^\pm $ and cups
  $ \mho^\pm $ will remain unchanged.

  From the uniqueness \cite{KhoroshkinTolstoy:91} of the universal
  $ R $ matrix for any quantum superalgebra the following relation
  holds (for appropriate normalisation):
  \be
    R^{-1} ( q )
    =
    R ( q^{-1} ),
  \ee
  which in turn leads to the relation
  \be
    \sigma^{-1} ( q )
    =
    P \sigma ( q^{-1} ) P.
  \ee
  Thus, up to a basis transformation, $ \sigma $ and $ \sigma^{-1} $ are
  interchangeable by the change of variable $ q \mapsto q^{-1} $. It
  then follows that the invariant for $ K^* $ is obtainable from that
  of $ K $ by the same change of variable, which leads to the
  following:

  \begin{proposition}
    If $ K $ is amphichiral then the invariant ${LG}_K$ is palindromic.%
    \footnote{%
      We intend ``palindromic'' to mean that the polynomial is
      invariant under the mapping $ q \mapsto q^{-1} $.
    }
    \proplabel{chiralmeanspalindromic}
  \end{proposition}

\vfill

\item[Inversion:]
  Again, if $ K^{-1} $ is the inverse of $ K $, then every arrow in
  $ K $ will have been replaced with an arrow in the opposite
  direction.  The tensor corresponding to $ K^{-1} $ will thus have the
  following changes: For the crossings,
  where $ X $ is either $ \sigma $ or $ \sigma^{-1} $, interchange
  $ X \iff X_d $ and $ X_l \iff X_r $;
  and for the caps and cups, interchange only the signs, i.e.
  $ \Omega^\pm \iff \Omega^\mp $
  and
  $ \mho^\pm \iff \mho^\mp $.

  This has the effect that the tensor representing $ K $ is replaced by
  the dual tensor acting on the dual space \cite{ReshetikhinTuraev:90}.
  Recalling that the tensors representing $ \( 1, 1 \) $ tangles act as
  scalar multiples of the identity on $ V $, then the dual tensor has
  exactly the same form, from which we conclude:

  \begin{proposition}
    A knot invariant derived from an irreducible representation of a
    quantum (super)algebra is unable to detect inversion.
    \proplabel{inversion}
  \end{proposition}
\end{description}


\subsection{Abstract Tensor Expressions for the Example Links}

We list the abstract tensors $ {\( T_K \)}^y_x $ that represent the
$ \( 1, 1 \) $-tangle (open diagram) forms of the example links.  In
each case, the indices $ x $ and $ y $ are the lower and upper loose
ends of the tangle in question.
The Links--Gould invariant is then formed by setting
$ x $ and $ y $ to be the same, i.e.
\be
  {LG}_K \( q, p = q^\alpha \)
  \defeq
  {\( T_K \)}^i_i
\ee
(no sum on $ i $), for any allowable index $ i $. We typically choose
$ i = 1 $.
Our invariant does not need to be \emph{writhe-normalised}, due to the
choice of normalisation of $ \sigma $ and the cap and cup matrices
$ \Omega^\pm $ and $ \mho^\pm $.  \figref{Loop} depicts removal of a
loop from a diagram.

\begin{figure}[htbp]
  \begin{center}
    \input{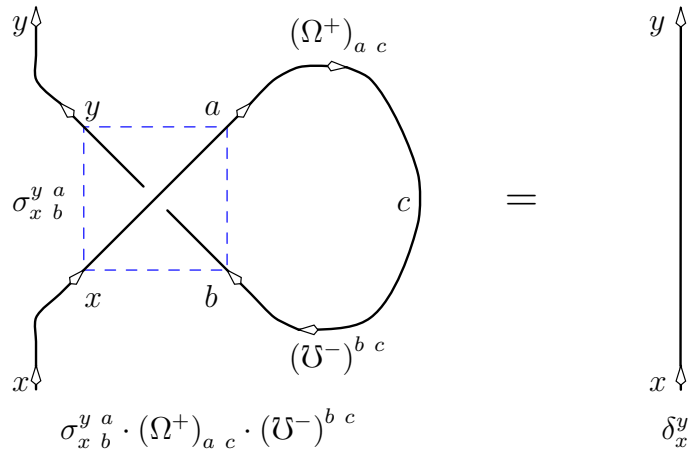}
    \caption{Removal of a Single, Positive Loop.}
    \figlabel{Loop}
  \end{center}
\end{figure}

Braid presentations for the example knots are taken
from \cite[pp~109-110]{Jones:85} and \cite[pp~381-386]{Jones:87}.

\begin{description}
\item[$ \mathbf{0_1} $ (Unknot):]
  A braid presentation is the trivial $ e \in B_1 $.
  As the unclosed tangle representing the Unknot is rather meaningless,
  we use simply $ {\( T_{0_1} \)}^y_x \defeq \delta^y_x $.

\item[$ \mathbf{2^2_1} $ (Hopf Link):]
  A braid presentation is $ {\sigma_1}^2 \in B_2 $.
  Diagrams pertaining to the Hopf Link and Trefoil are found in
  \figref{HopfLinkTrefoil}.
  \be
    {( T_{2^2_1} )}^y_x
    \defeq
    {\( \sigma^2 \)}^{y~a}_{x~b}
    \cdot
    {\( \Omega^+ \)}_{a~c}
    \cdot
    {\( \mho^- \)}^{b~c}.
  \ee

\item[$ \mathbf{3_1} $ (Trefoil):]
  A braid presentation is $ {\sigma_1}^3 \in B_2 $.
  This knot has also been
  called the \emph{overhand knot} (as that is how it is tied) and the
  \emph{cloverleaf knot} \cite[pp~3-4]{CrowellFox:77}.
  \be
    {\( T_{3_1} \)}^y_x
    \defeq
    {\( \sigma^3 \)}^{y~a}_{x~b}
    \cdot
    {\( \Omega^+ \)}_{a~c}
    \cdot
    {\( \mho^- \)}^{b~c}.
  \ee

  \begin{figure}[htbp]
    \begin{center}
      \input{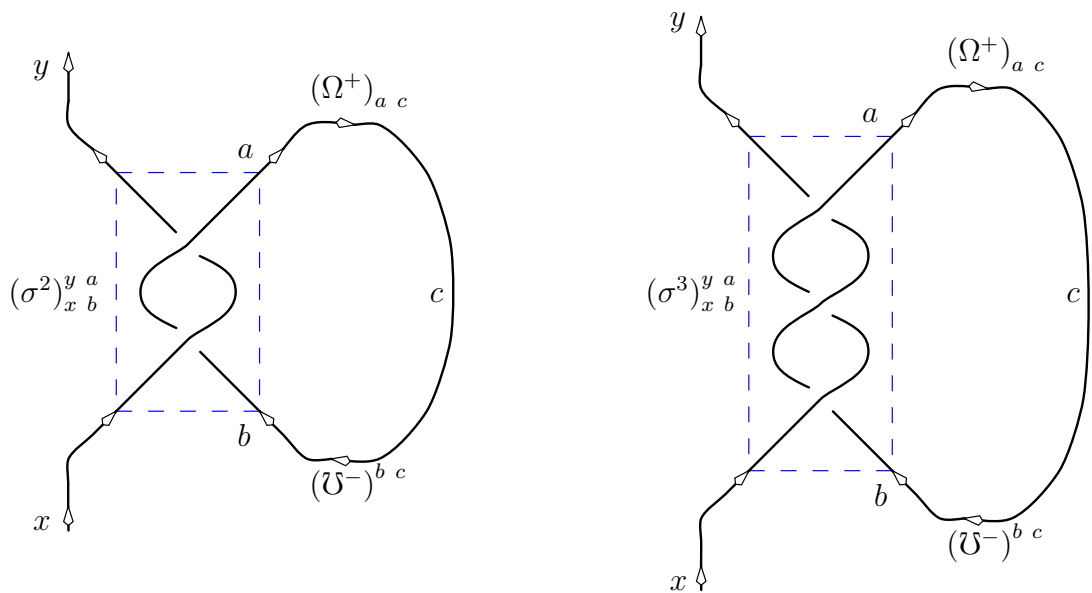}
      \caption{Tangle form of $ 2^2_1 $ (the Hopf Link) and $ 3_1 $
         (the positive Trefoil).
      }
      \figlabel{HopfLinkTrefoil}
    \end{center}
  \end{figure}

\pagebreak

\item[$ \mathbf{4_1} $ (Figure Eight):]
  A braid presentation is
  $ {\( \sigma_1 \sigma_2^{-1} \)}^2 \in B_3 $, and
  a diagram is found in \figref{FigureEight}. This knot has also been
  called the \emph{Four-Knot} (as it is the only $ 4 $ crossing knot)
  and \emph{Listing's Knot} \cite[p~4]{CrowellFox:77}.
  \be
    {\( T_{4_1} \)}^y_x
    \defeq
    {\( \sigma^{-1}_l \sigma^{-1}_r \)}^{y~b}_{a~c}
    \cdot
    {\( \sigma_r \)}^{c~e}_{d~f}
    \cdot
    {\sigma}^{a~d}_{x~g}
    \cdot
    {\( \Omega^- \)}_{b~e}
    \cdot
    {\( \mho^- \)}^{g~f}.
  \ee

  \begin{figure}[ht]
    \begin{center}
      \input{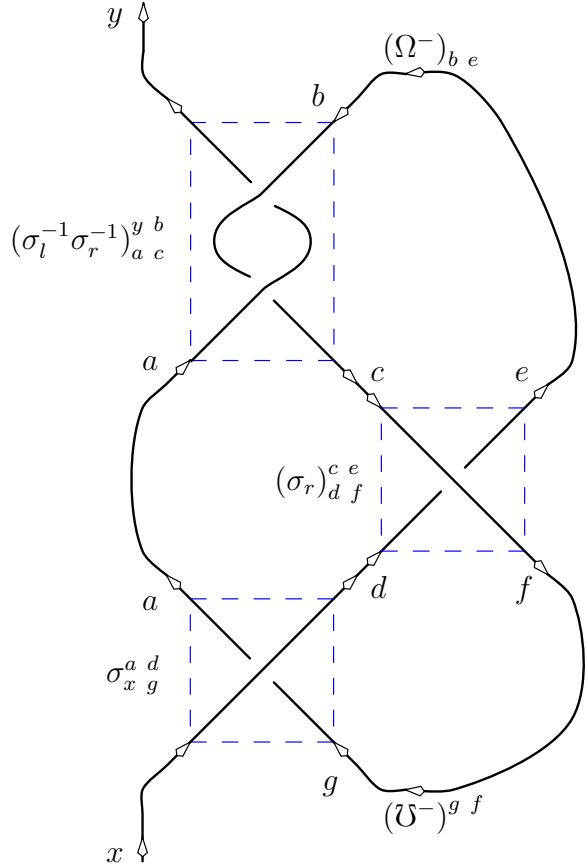}
      \caption{Tangle form of $ 4_1 $ (the Figure Eight Knot).
      }
      \figlabel{FigureEight}
    \end{center}
  \end{figure}

\item[$ \mathbf{5^2_1} $ (Whitehead Link):]
  A braid presentation is
  $ {\( \sigma_1 \sigma_2^{-1} \)}^2 \sigma_2^{-1} \in B_3 $ and a
  diagram is found in \figref{WhiteheadLink}. (This link is named after
  the topologist J H C Whitehead, not the logician Alfred North
  Whitehead \cite[p~200]{Kauffman:88}.)

  Firstly, we define a temporary tensor to reduce computation:
  \be
    {\( W \)}^{c~i}_{x~d}
    \defeq
    {\( \sigma^{-2} \)}^{c~e}_{x~f}
    \cdot
    {( \sigma_d^2 )}^{g~i}_{h~d}
    \cdot
    {\( \Omega^+ \)}_{e~g}
    \cdot
    {\( \mho^- \)}^{f~h},
  \ee
  where we have written $ \sigma^{-2} \defeq {\( \sigma^{-1} \)}^2 $.
  With this, we have:
  \be
    {( T_{5^2_1} )}^y_x
    \defeq
    {\( W \)}^{c~i}_{x~d}
    \cdot
    {\( \sigma_r \sigma_l \)}^{a~y}_{i~b}
    \cdot
    {\( \Omega^+ \)}_{c~a}
    \cdot
    {\( \mho^+ \)}^{d~b}.
  \ee

  \begin{figure}[htbp]
    \begin{center}
      \input{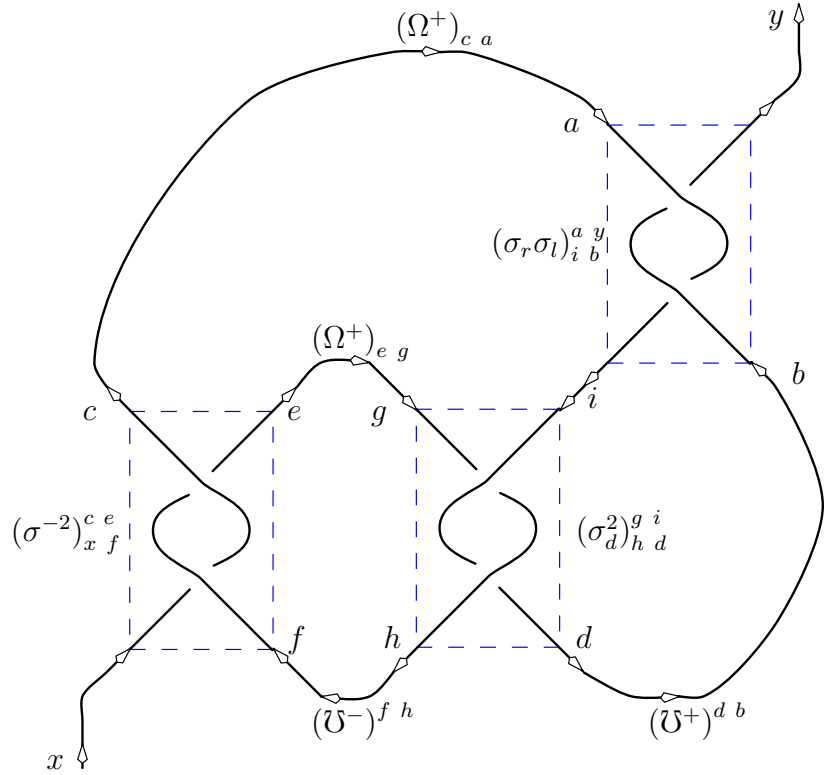}
      \caption{Tangle form of $ 5^2_1 $ (the Whitehead Link).
      }
      \figlabel{WhiteheadLink}
    \end{center}
  \end{figure}

\item[$ \mathbf{8_{17}} $:]
  A braid presentation is
  $
    {\( \sigma_1^{-1} \sigma_2 \)}^2 \sigma_2^2 \sigma_1^{-2} \sigma_2
    \in B_3
  $, and a diagram is found in \figref{EightSeventeen}.
  Again, we define some temporary tensors to reduce computation:
  \be
    {\( EA \)}^{y~c~e}_{b~d~f}
    & \defeq &
    {\( \sigma^{-2} \)}^{c~e}_{g~f}
    \cdot
    {\( \sigma^2 \)}^{y~g}_{b~d}
    \\
    {\( EB \)}^{b~d~f}_{x~i~j}
    & \defeq &
    {{\( \sigma^{-1} \)}}^{d~f}_{k~l}
    \cdot
    {\sigma}^{b~k}_{m~n}
    \cdot
    {{\( \sigma^{-1} \)}}^{n~l}_{o~j}
    \cdot
    {\sigma}^{m~o}_{x~i}.
  \ee
  With these, we have:
  \be
    {\( T_{8_{17}} \)}^y_x
    \defeq
    {\( EA \)}^{y~c~e}_{b~d~f}
    \cdot
    {\( EB \)}^{b~d~f}_{x~i~j}
    \cdot
    {\( \Omega^+ \)}_{c~r}
    \cdot
    {\( \mho^- \)}^{i~r}
    \cdot
    {\( \Omega^+ \)}_{e~q}
    \cdot
    {\( \mho^- \)}^{j~q}.
  \ee
  To reduce computation, we may define even more auxiliary tensors:
  \be
    {\( EB \)}^{b~d~f}_{x~i~j}
    =
    {\( EC \)}^{b~d~f}_{m~n~l}
    \cdot
    {\( ED \)}^{m~n~l}_{x~i~j},
  \ee
  where:
  \be
    {\( EC \)}^{b~d~f}_{m~n~l}
    & \defeq &
    {{\( \sigma^{-1} \)}}^{d~f}_{k~l}
    \cdot
    {\sigma}^{b~k}_{m~n}
    \\
    {\( ED \)}^{m~n~l}_{x~i~j}
    & \defeq &
    {{\( \sigma^{-1} \)}}^{n~l}_{o~j}
    \cdot
    {\sigma}^{m~o}_{x~i}.
  \ee

  \begin{figure}[htbp]
    \begin{center}
      \input{EightSeventeen.pstex_t}
      \caption{Tangle form of $ 8_{17} $.}
      \figlabel{EightSeventeen}
    \end{center}
  \end{figure}

\pagebreak

\item[$ \mathbf{9_{42}} $:]
  A braid presentation is
  $
    \sigma_1^3 \sigma_3
    \sigma_2^{-1} \sigma_3
    \sigma_1^{-2} \sigma_2^{-1}
    \in B_4
  $,
  and a diagram is found in \figref{NineFortytwo}.
  Again, we define a temporary tensor to reduce computation:
  \be
    {\( N \)}^{a~y}_{b~h}
    & \defeq &
    {( \sigma_d^2 )}^{a~c}_{b~d}
    \cdot
    {\( \sigma^{-3} \)}^{e~y}_{f~h}
    \cdot
    {\( \Omega^- \)}_{c~e}
    \cdot
    {\( \mho^+ \)}^{d~f}.
    \\
    {\( T_{9_{42}} \)}^y_x
    & \defeq &
    {\( N \)}^{a~y}_{b~h}
    \cdot
    {\( \sigma^{-1}_d \sigma^{-1} \)}^{b~h}_{i~j}
    \cdot
    {\( \sigma \sigma_d \)}^{k~i}_{x~m}
    \cdot
    {\( \mho^+ \)}^{m~j}
    \cdot
    {\( \Omega^+ \)}_{k~a}.
  \ee

\item[$ \mathbf{10_{48}} $:]
  A braid presentation is
  $
    \sigma_1^{-2} \sigma_2^{4} \sigma_1^{-3} \sigma_2
    \in B_3
  $,
  and a diagram is found in \figref{TenFortyEight}.
  \be
    {\( T_{10_{48}} \)}^y_x
    \defeq
    {{\( \sigma^{-2} \)}}^{a~y}_{b~f}
    \cdot
    {\( \sigma^{4} \)}^{f~g}_{d~h}
    \cdot
    {\( \sigma^{-3} \)}^{b~d}_{c~e}
    \cdot
    {\( \sigma \)}^{e~h}_{x~i}
    \cdot
    \\
    \qquad \qquad
    {\( \Omega^- \)}_{j~a}
    \cdot
    {\( \mho^+ \)}^{j~c}
    \cdot
    {\( \Omega^+ \)}_{g~k}
    \cdot
    {\( \mho^- \)}^{i~k}.
  \ee

  \begin{figure}[htbp]
    \begin{center}
      \input{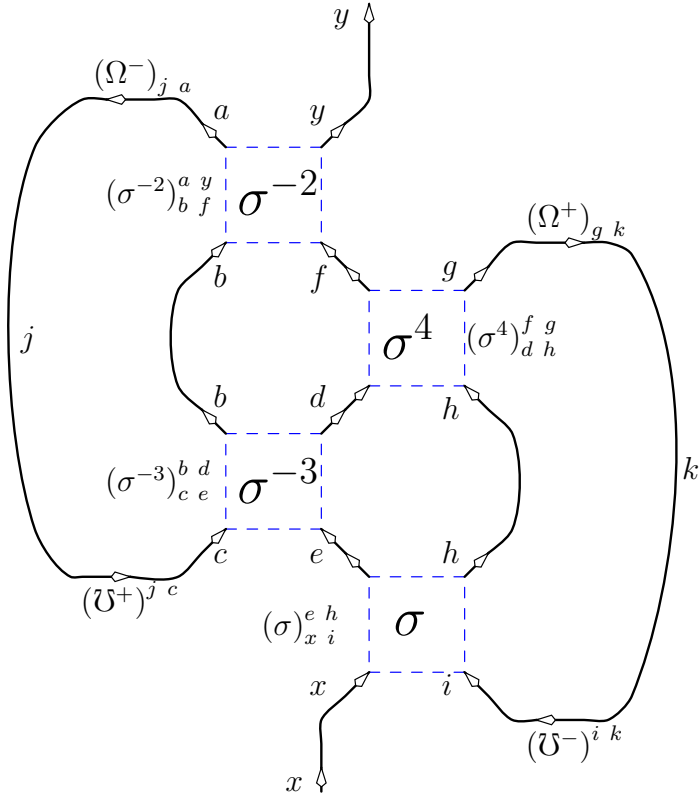}
      \caption{Tangle form of $10_{48}$.}
      \figlabel{TenFortyEight}
    \end{center}
  \end{figure}

  \begin{figure}[ht]
    \begin{center}
      \input{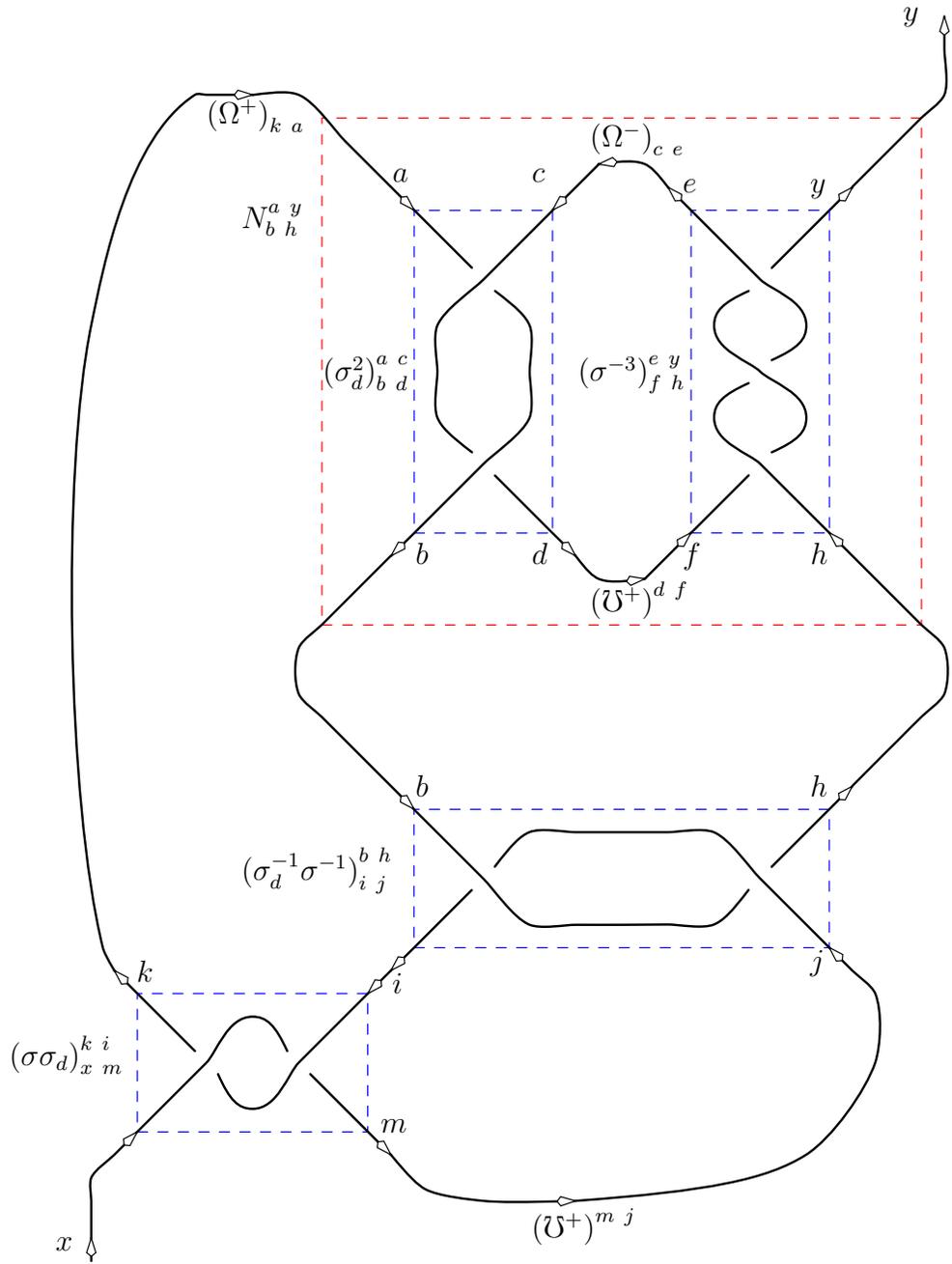}
      \caption{Tangle form of $ 9_{42} $.}
      \figlabel{NineFortytwo}
    \end{center}
  \end{figure}
\end{description}

\clearpage


\section{The Links--Gould Tangle Invariant}
\seclabel{lglink}

\subsection{Crossing Matrices $ \sigma $ and $ \sigma^{-1} $}

From the results of \secref{Rmatrixdef}, we have the matrices
$ \sigma $ and $ \sigma^{-1} $,
using the substitution $ p \defeq q^\alpha $:
\be
  \sigma
  =
  \m{[}{@{}*{3}{*{3}{c@{}@{}}c|}*{3}{c@{}@{}}c@{}}
    {\scs p^{-2}} &. &. &. &. &. &. &. &. &. &. &. &. &. &. &. \\
    . &. &. &. &{\scs p^{-1}} &. &. &. &. &. &. &. &. &. &. &. \\
    . &. &. &. &. &. &. &. &{\scs p^{-1}} &. &. &. &. &. &. &. \\
    . &. &. &. &. &. &. &. &. &. &. &. &{\scs1} &. &. &. \\
    \hline
    . &{\scs p^{-1}} &. &. &{\scs p^{-2}-1} &. &. &. &. &. &. &. &. &. &. &. \\
    . &. &. &. &. &{\scs -1} &. &. &. &. &. &. &. &. &. &. \\
    . &. &. &. &. &. &{\scs q^2-1} &. &. &{\scs -q} &. &. &{\scs - q Y} &. &. &. \\
    . &. &. &. &. &. &. &. &. &. &. &. &. &{\scs p q} &. &. \\
    \hline
    . &. &{\scs p^{-1}} &. &. &. &. &. &{\scs p^{-2}-1} &. &. &. &. &. &. &. \\
    . &. &. &. &. &. &{\scs -q} &. &. &. &. &. &{\scs Y} &. &. &. \\
    . &. &. &. &. &. &. &. &. &. &{\scs -1} &. &. &. &. &. \\
    . &. &. &. &. &. &. &. &. &. &. &. &. &. &{\scs p q} &. \\
    \hline
    . &. &. &{\scs 1} &. &. &{\scs - q Y} &. &. &{\scs Y} &. &. &{\scs Y^2} &. &. &. \\
    . &. &. &. &. &. &. &{\scs p q} &. &. &. &. &. &{\scs p^2 q^2 -1} &. &. \\
    . &. &. &. &. &. &. &. &. &. &. &{\scs p q} &. &. &{\scs p^2 q^2 -1} &. \\
    . &. &. &. &. &. &. &. &. &. &. &. &. &. &. &{\scs p^2 q^2}
  \me,
\ee
\be
  \! \! \! \! \! \! \! \! \! \! \! \! \!
  \! \! \!
  \sigma^{-1}
  =
  \m{[}{@{}*{3}{*{3}{c@{}@{}}c|}*{3}{c@{}@{}}c@{}}
    {\scs p^2} &. &. &. &. &. &. &. &. &. &. &. &. &. &. &. \\
    . &{\scs p^2 -1} &. &. &{\scs p} &. &. &. &. &. &. &. &. &. &. &. \\
    . &. &{\scs p^2 -1} &. &. &. &. &. &{\scs p} &. &. &. &. &. &. &. \\
    . &. &. &{\scs Y^2 q^{-2}} &. &. &{\scs Y q^{-1}} &. &. &{\scs -Y q^{-2}} &. &. & {\scs 1} &. &. &. \\
    \hline
    . &{\scs p} &. &. &. &. &. &. &. &. &. &. &. &. &. &. \\
    . &. &. &. &. &{\scs -1} &. &. &. &. &. &. &. &. &. &. \\
    . &. &. &{\scs Y q^{-1}} &. &. &. &. &. &{\scs -q^{-1}} &. &. &. &. &. &. \\
    . &. &. &. &. &. &. &{\scs p^{-2} q^{-2} - 1} &. &. &. &. &. &{\scs p^{-1} q^{-1}} &. & . \\
    \hline
    . &. &{\scs p} &. &. &. &. &. &. &. &. &. &. &. &. &. \\
    . &. &. &{\scs -Y q^{-2}} &. &. &{\scs -q^{-1}} &. &. &{\scs q^{-2} - 1} &. &. &. &. &. & . \\
    . &. &. &. &. &. &. &. &. &. &{\scs -1} &. &. &. &. &. \\
    . &. &. &. &. &. &. &. &. &. &. &{\scs p^{-2} q^{-2} - 1} &. &. &{\scs p^{-1} q^{-1}} & . \\
    \hline
    . &. &. &{\scs 1} &. &. &. &. &. &. &. &. &. &. &. &. \\
    . &. &. &. &. &. &. &{\scs p^{-1} q^{-1}} &. &. &. &. &. &. &. &. \\
    . &. &. &. &. &. &. &. &. &. &. &{\scs p^{-1} q^{-1}} &. &. &. &. \\
    . &. &. &. &. &. &. &. &. &. &. &. &. &. &. &{\scs p^{-2} q^{-2}}
  \me,
\ee
where $ Y = {\( p^{-2} - q^2 + p^2 q^2 - 1 \)}^{1/2} $.


\subsection{Caps and Cups $ \Omega^\pm $ and $ \mho^\pm $}

Where $ \mho^\pm = {\( \Omega^\pm \)}^{-1} $, we will use:
\be
  \Omega^-
  \eq
  \m{[}{c@{}c@{}c@{}c}
    {\scs q^{- 2 \alpha}} & . & . & . \\
    . & {\scs - q^{- 2 \( \alpha + 1 \)}} & . & . \\
    . & . & {\scs - q^{- 2 \alpha}} & . \\
    . & . & . & {\scs q^{- 2 \( \alpha + 1 \)}}
  \me
  \equiv
  \m{[}{cccc}
    {\scs p^{-2}} & . & . & . \\
    . & {\scs - p^{-2} q^{-2}} & . & . \\
    . & . & {\scs - p^{-2}} & . \\
    . & . & . & {\scs p^{-2} q^{-2}}
  \me,
  \\
  \qquad
  \mho^-
  \eq
  \m{[}{c@{}c@{}c@{}c}
    {\scs q^{2 \alpha}} & . & . & . \\
    . & {\scs - q^{2 \( \alpha + 1 \)}} & . & . \\
    . & . & {\scs - q^{2 \alpha}} & . \\
    . & . & . & {\scs q^{2 \( \alpha + 1 \)}}
  \me
  \equiv
  \m{[}{cccc}
    {\scs p^2} &               . &        . & . \\
    .      & {\scs - p^2 q^2} &        . & . \\
    .      &               . & {\scs - p^2} & . \\
    .      &               . &        . & {\scs p^2 q^2}
  \me,
  \\
  \Omega^+
  \eq
  \mho^+
  =
  I_4.
\ee

The choices for $ \Omega^\pm $ and $ \mho^\pm $ are not unique.

\begin{itemize}
\item
  $ \Omega^+ $ and $ \mho^+ $ may be chosen from
  consistency considerations in \figref{XLandXR}. The simple choices:
  \be
    {\( \Omega^+ \)}_{a~b}
    =
    {\( \mho^+ \)}^{a~b}
    =
    \delta_{a~b}
  \ee
  (i.e. $ \Omega^+ = \mho^+ = I_4 $),
  ensure that the definition \eqref{XRdef}, i.e.
  \be
    {\( X_r \)}^{a~c}_{b~d}
    \defeq
    {X}^{c~g}_{f~b}
    \cdot
    {\( \mho^+ \)}^{a~f}
    \cdot
    {\( \Omega^+ \)}_{g~d}
  \ee
  (where $ X $ is either $ \sigma $ or $ \sigma^{-1} $),
  simplifies to the elegant form:
  \be
    {\( X_r \)}^{a~c}_{b~d}
    =
    {X}^{c~d}_{a~b}.
  \ee

\item
  For the choice of $ \Omega^- $ and $ \mho^- $,
  we invoke the following result from
  \cite[Lemma~2,~p~354]{LinksGouldZhang:93}
  (see also \cite{LinksGould:92b}):
  \be
    \( I \otimes \mathrm{str} \)
    [ ( I \otimes q^{-2 h_\rho} ) \sigma ]
    =
    k I,
  \ee
  for some constant $ k $ depending on the normalisation of $ \sigma $.
  Note that $ \mathrm{str} $ denotes the supertrace, and that in this
  case:
  \be
    \pi ( q^{-2 h_\rho} )
    =
    \m{[}{cccc}
      q^{- 2 \alpha} &              . &              . &              . \\
               . & q^{- 2 \alpha - 2} &              . &              . \\
               . &              . & q^{- 2 \alpha}     &              . \\
               . &              . &              . & q^{- 2 \alpha - 2}
    \me.
  \ee

  \pagebreak

  From \figref{Loop}, we require:
  \be
    \sigma^{y a}_{x b}
    \cdot
    {\( \Omega^+ \)}_{a c}
    \cdot
    {\( \mho^- \)}^{b c}
    =
    \delta^y_x,
  \ee
  which, along with the condition:
  \be
    {\( \Omega^- \)}_{a b}
    \cdot
    {\( \mho^- \)}^{b c}
    =
    \delta^c_a,
  \ee
  imposes the choice:
  \be
    {\( \mho^- \)}^{b c}
    \eq
    {\( - \)}^{\[ b \]}
    \pi {( q^{-2 h_\rho} )}^b_c,
    \\
    {\( \Omega^- \)}_{b c}
    \eq
    {\( - \)}^{\[ b \]}
    \pi {( q^{2 h_\rho} )}^b_c.
  \ee
\end{itemize}

For other references on the construction of the cap and cup matrices,
see the papers by Reshetikhin and Turaev \cite{ReshetikhinTuraev:90},
and particularly Zhang \cite{Zhang:95} for the superalgebra case.


\subsection{Results}

Some evaluations of the invariant are presented in
\tabref{gl21oriented}.  This $ U_q \[ gl \( 2 | 1 \) \] $ oriented
invariant \emph{is} an invariant of ambient isotopy.

\begin{table}[ht]
  \centering
  \begin{tabular}{||c|l||}
    \hline\hline
    & \\[-3mm]
    $ K $ &
    \multicolumn{1}{c||}{$ {LG}_K \( q, p \) $}
    \\[1mm]
    \hline\hline
    & \\[-2mm]
    $ 0_1 $ &
    $ 1 $
    \\[2mm]
    \hline
    & \\[-2mm]
    $ 2^2_1 $ &
    $ -1 + p^{-2} - q^2 + p^2 q^2 $
    \\[2mm]
    \hline
    & \\[-2mm]
    $ 3_1 $ &
    $
      1 + p^{-4} - p^{-2} + 2 q^2 - p^{-2} q^2 - p^2 q^2
      - p^2 q^4 + p^4 q^4
    $
    \\[2mm]
    \hline
    & \\[-2mm]
    $ 4_1 $ &
    $
        7
      +   \( p^{-4} q^{-2} + p^4 q^2 \)
      - 3 \( p^{-2}        + p^2     \)
      - 3 \( p^{-2} q^{-2} + p^2 q^2 \)
      + 2 \(        q^{-2} +     q^2 \)
    $
    \\[2mm]
    \hline
    & \\[-2mm]
    $ 5^2_1 $ &
    $
      \begin{array}{l}
        - 10
        +   p^{-6} q^{-2}
        - 3 p^{-4}
        - 3 p^{-4} q^{-2}
        + 4 p^{-2} q^{-2}
        + 9 p^{-2}
        - 2 q^{-2}
        \\
        \quad
        - 8        q^2
        + 2 p^{-2} q^2
        + 9 p^2    q^2
        + 4 p^2
        + 2 p^2    q^4
        - 3 p^4    q^2
        - 3 p^4    q^4
        +   p^6    q^4
      \end{array}
    $
    \\[2mm]
    \hline\hline
    & \\[-2mm]
    $ 8_{17} $ &
    see \secref{EightSeventeen}
    \\[2mm]
    \hline
    & \\[-2mm]
    $ 9_{42}, 10_{48} $ &
    see \secref{NineFortyTwoandTenFortyEight}
    \\[2mm]
    \hline\hline
  \end{tabular}
  \caption{%
    The Links--Gould $ U_q \[ gl \( 2 | 1 \) \] $ oriented
    polynomial invariant $ {LG}_K \( q, p \) $,
    evaluated using the open diagram form of various links $ K $.}
  \tablabel{gl21oriented}
\end{table}


\subsection{Behaviour of the Invariant}

Fix a knot $ K $, and denote by $ K^* $ the reflection of $ K $ and by
$ K^{-1} $ the inverse of $ K $.  From the polynomial for $ K $, we may
immediately write down the polynomials for $ K^* $ and $ K^{-1} $.  For
the reflection, we have:
\bse
  {LG}_{K^*} \( q, p \)
  =
  {LG}_K \( q^{-1}, p^{-1} \).
  \eqlabel{Utricularia.1}
\ese
For the inverse, we have:
\be
  {LG}_{K^{-1}} \( q, p \)
  =
  {LG}_K \( q, q^{-1} p^{-1} \).
\ee
(this follows from $ \alpha \mapsto - \( \alpha + 1 \) $).

\begin{description}
\item[Chirality:]
  As we have:
  \be
    K = K^*
    \quad \to \quad
    {LG}_K \( q, p \)
    =
    {LG}_{K^*} \( q, p \),
  \ee
  then we have, conversely, that:
  \bse
    {LG}_K \( q, p \)
    \neq
    {LG}_{K^*} \( q, p \)
    \quad \to \quad
    K \neq K^*,
    \eqlabel{Utricularia.3}
  \ese
  i.e. if the polynomials corresponding to $ K $ and $ K^* $ are
  distinct, then $ K $ must be chiral. Using the identity
  \eqref{Utricularia.1}, the test of \eqref{Utricularia.3} becomes:
  \be
    {LG}_K \( q, p \)
    \neq
    {LG}_K \( q^{-1}, p^{-1} \)
    \quad \to \quad
    K \neq K^*,
  \ee
  i.e. if $ {LG}_K \( q, p \) $
  is \emph{not} palindromic, then $ K $ is chiral.
\item[Invertibility:]
  We make the observation that the representation of
  $ U_q \[ gl \( 2 | 1 \) \] $ acting on the dual module $ V^* $
  is given by the replacement $ \alpha \mapsto - \( \alpha + 1 \) $
  (with an appropriate redefinition of the Cartan elements). Thus
  for a given $ \( 1, 1 \) $ tangle $ K $, with invariant
  $ {LG}_K \( q, p \) $, the tangle invariant
  $ {LG}_{K^{-1}} $ of its inverse $ K^{-1} $ is obtained as
  $ {LG}_{K^{-1}} \( q, p \) = {LG}_K \( q, q^{-1} p^{-1} \) $.
  However, in view of \propref{inversion}, such an invariant is unable
  to detect inversion.
\end{description}

We summarise these results in a proposition:

\begin{proposition}
  Let $ {LG}_K \( q, p \) $ be the Links--Gould polynomial invariant for
  the knot $ K $.
  \begin{itemize}
  \item
    If $ {LG}_K \( q, p \) $ is \emph{not}
    invariant under the transformation $ q \mapsto q^{-1} $
    (which implies $ p \mapsto p^{-1} $), then
    $ LG $ detects the chirality of $ K $.
  \item
    $ {LG}_K \( q, p \) $ enjoys the symmetry property:
    \bse
      {LG}_K \( q, p \) = {LG}_K \( q, q^{-1} p^{-1} \).
      \eqlabel{Melanotaenia}
    \ese
  \end{itemize}
  \proplabel{Rhadinocentrus}
\end{proposition}


\subsection{The Chirality of $9_{42}$ and $10_{48}$}
\seclabel{NineFortyTwoandTenFortyEight}

The polynomials for $9_{42}$ and $10_{48}$ are:
\begin{eqnarray*}
  & &
  \! \! \! \! \! \! \! \! \! \! \! \! \! \! \! \!
  \! \! \! \! \! \! \! \!
  {LG}_{9_{42}} \( q, p \)
  =
  \\
  & &
  \! \! \! \! \! \! \! \! \! \! \! \! \! \! \! \!
  \begin{array}{l}
      3
    +   p^{-8} q^{-6}
    - 2 p^{-6} q^{-6}
    - 2 p^{-6} q^{-4}
    +   p^{-4} q^{-6}
    + 3 p^{-4} q^{-4}
    +   p^{-4} q^{-2}
    +   p^{-4}
    -   p^{-2} q^{-4}
    \\
    -   p^{-2} q^{-2}
    - 3 p^{-2}
    - 3 p^{-2} q^2
    + 6        q^2
    + 2        q^4
    -   p^2    q^{-2}
    -   p^2
    - 3 p^2    q^2
    - 3 p^2    q^4
    +   p^4    q^{-2}
    \\
    + 3 p^4
    +   p^4    q^2
    +   p^4    q^4
    - 2 p^6
    - 2 p^6    q^2
    +   p^8    q^2
  \end{array}
  \\
  & &
  \! \! \! \! \! \! \! \! \! \! \! \! \! \! \! \!
  \! \! \! \! \! \! \! \!
  {LG}_{10_{48}} \( q, p \)
  =
  \\
  & &
  \! \! \! \! \! \! \! \! \! \! \! \! \! \! \! \!
  \begin{array}{l}
      165
    + 5   p^{-8}
    - 25  p^{-6}
    + 68  p^{-4}
    - 129 p^{-2}
    - 132 p^{2}
    + 67  p^{4}
    - 22  p^{6}
    + 4   p^{8}
    +     p^{-16} q^{-8}
    \\
    - 3   p^{-14} q^{-8}
    + 4   p^{-12} q^{-8}
    - 4   p^{-10} q^{-8}
    + 4   p^{-8}  q^{-8}
    - 2   p^{-6}  q^{-8}
    - 3   p^{-14} q^{-6}
    + 12  p^{-12} q^{-6}
    \\
    - 21  p^{-10} q^{-6}
    + 24  p^{-8}  q^{-6}
    - 22  p^{-6}  q^{-6}
    + 13  p^{-4}  q^{-6}
    - 3   p^{-2}  q^{-6}
    + 16          q^{-4}
    + 5   p^{-12} q^{-4}
    \\
    - 23  p^{-10} q^{-4}
    + 50  p^{-8}  q^{-4}
    - 69  p^{-6}  q^{-4}
    + 67  p^{-4}  q^{-4}
    - 43  p^{-2}  q^{-4}
    - 3   p^{2}   q^{-4}
    + 94          q^{-2}
    \\
    - 6   p^{-10} q^{-2}
    + 29  p^{-8}  q^{-2}
    - 72  p^{-6}  q^{-2}
    + 119 p^{-4}  q^{-2}
    - 132 p^{-2}  q^{-2}
    - 43  p^{2}   q^{-2}
    + 13  p^{4}   q^{-2}
    \\
    - 2   p^{6}   q^{-2}
    + 88          q^{2}
    - 2   p^{-6}  q^{2}
    + 12  p^{-4}  q^{2}
    - 39  p^{-2}  q^{2}
    - 129 p^{2}   q^{2}
    + 119 p^{4}   q^{2}
    - 69  p^{6}   q^{2}
    \\
    + 24  p^{8}   q^{2}
    - 4   p^{10}  q^{2}
    + 12          q^{4}
    - 2   p^{-2}  q^{4}
    - 39  p^{2}   q^{4}
    + 68  p^{4}   q^{4}
    - 72  p^{6}   q^{4}
    + 50  p^{8}   q^{4}
    - 21  p^{10}  q^{4}
    \\
    + 4   p^{12}  q^{4}
    - 2   p^{2}   q^{6}
    + 12  p^{4}   q^{6}
    - 25  p^{6}   q^{6}
    + 29  p^{8}   q^{6}
    - 23  p^{10}  q^{6}
    + 12  p^{12}  q^{6}
    - 3   p^{14}  q^{6}
    - 2   p^{6}   q^{8}
    \\
    + 5   p^{8}   q^{8}
    - 6   p^{10}  q^{8}
    + 5   p^{12}  q^{8}
    - 3   p^{14}  q^{8}
    +     p^{16}  q^{8}.
  \end{array}
\end{eqnarray*}
Neither of these polynomials are palindromic, hence $LG$ \emph{does}
distinguish the chirality of these knots.


\subsection{The Noninvertibility of $ 8_{17} $ is not Detected}
\seclabel{EightSeventeen}

Recall that $ 8_{17} $ is the smallest noninvertible knot.
We find its polynomial invariant to be given by:
\be
  & &
  \! \! \! \! \! \! \! \! \! \! \! \! \! \! \! \!
  \! \! \! \! \! \! \! \!
  {LG}_{8_{17}} \( q, p \)
  =
  \\
  & &
  \! \! \! \! \! \! \! \! \! \! \! \! \! \! \! \!
  \begin{array}{l}
      139
    +     \( p^{-12} q^{-6} + p^{12} q^6    \)
    -   4 \( p^{-10} q^{-6} + p^{10} q^6    \)
    -   4 \( p^{-10} q^{-4} + p^{10} q^4    \)
    \\
    +   7 \( p^{-8}  q^{-6} + p^8    q^6    \)
    +  18 \( p^{-8}  q^{-4} + p^8    q^4    \)
    +   7 \( p^{-8}  q^{-2} + p^8    q^2    \)
    -   7 \( p^{-6}  q^{-6} + p^6    q^6    \)
    \\
    -  36 \( p^{-6}  q^{-4} + p^6    q^4    \)
    -  36 \( p^{-6}  q^{-2} + p^6    q^2    \)
    -   7 \( p^{-6}         + p^6           \)
    +   3 \( p^{-4}  q^{-6} + p^4    q^6    \)
    \\
    +  40 \( p^{-4}  q^{-4} + p^4    q^4    \)
    +  82 \( p^{-4}  q^{-2} + p^4    q^2    \)
    +  40 \( p^{-4}         + p^4           \)
    +   3 \( p^{-4}  q^2    + p^4    q^{-2} \)
    \\
    -  22 \( p^{-2}  q^{-4} + p^2    q^4    \)
    - 102 \( p^{-2}  q^{-2} + p^2    q^2    \)
    - 102 \( p^{-2}         + p^2           \)
    -  22 \( p^{-2}  q^2    + p^2    q^{-2} \)
    \\
    +   4 \( q^{-4}         + q^4           \)
    +  68 \( q^{-2}         + q^2           \).
  \end{array}
\ee
As $ 8_{17} $ is amphichiral, the polynomial invariant is palindromic,
as predicted by \propref{Rhadinocentrus}.
Furthermore, we may observe the invariance:
\be
  {LG}_{8_{17}} \( q, p \)
  =
  {LG}_{8_{17}} \( q, q^{-1} p^{-1} \),
\ee
which is consistent with our assertion that our polynomial invariant
cannot detect the noninvertibility of \emph{any} knot. More
experiments to illustrate this claim are supplied in
\secref{trotterpretzel}.


\subsection{A Class of Noninvertible Pretzel Knots}
\seclabel{trotterpretzel}

A class of noninvertible knots has been presented by Trotter
\cite{Trotter:64}; this class provides an easily-programmable
set of examples to see if a knot invariant detects noninvertibility.
Trotter is of the opinion that the knots are chiral.  These pretzel
knots were in fact the \emph{first} noninvertible knots to be described
\cite[p~25]{Livingston:93}.

The structure of the knots $ \( p, q, r \) $ in this family is depicted
by its simplest example in \figref{TrotterPretzel}. Note that $ p, q $,
and $ r $ must all be distinct, odd, and greater than $ 1 $.  In
\figref{TrotterPretzel}, the notation $ X_{rlr}^N $ refers to the
$ 2 $-braid of $ N $ crossings formed by the placing of $ X_r $ atop
$ X_l $, with $ X_r $ as the top and bottom crossings, for $ X $ being
either $ \sigma $ or $ \sigma^{-1} $. The recursive definition for such
knots is provided in \figref{XRLR}.

\begin{figure}[htbp]
  \begin{center}
    \input{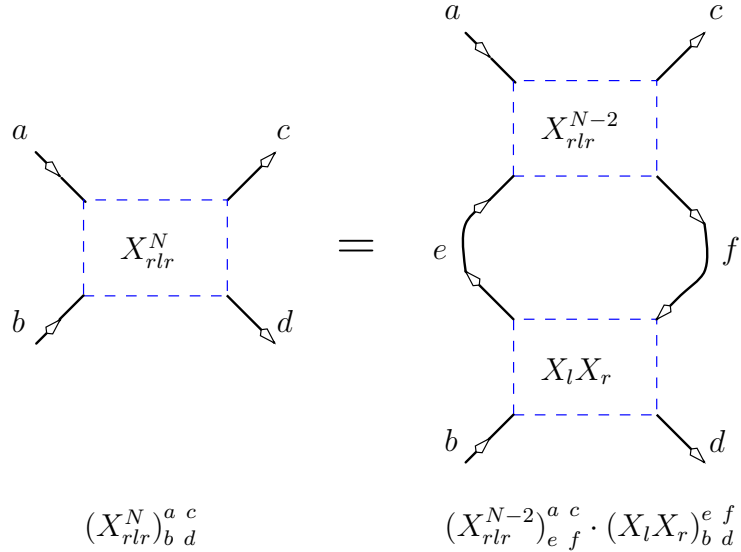}
    \caption{%
       Recursive definition of the towers $ X_{rlr} $ used in the
       evaluation of the Links--Gould link invariant for the Trotter
       pretzel knots; $X$ is either $ \sigma $ or $ \sigma^{-1} $.
       The minimum is the case $ N = 1 $, which
       corresponds to $ X_r $, i.e.  $ X_{rlr}^1 \defeq X_r $.
       A parallel definition of $ X_{lrl} $ might be given.
    }
    \figlabel{XRLR}
  \end{center}
\end{figure}

\begin{figure}[htbp]
  \begin{center}
    \input{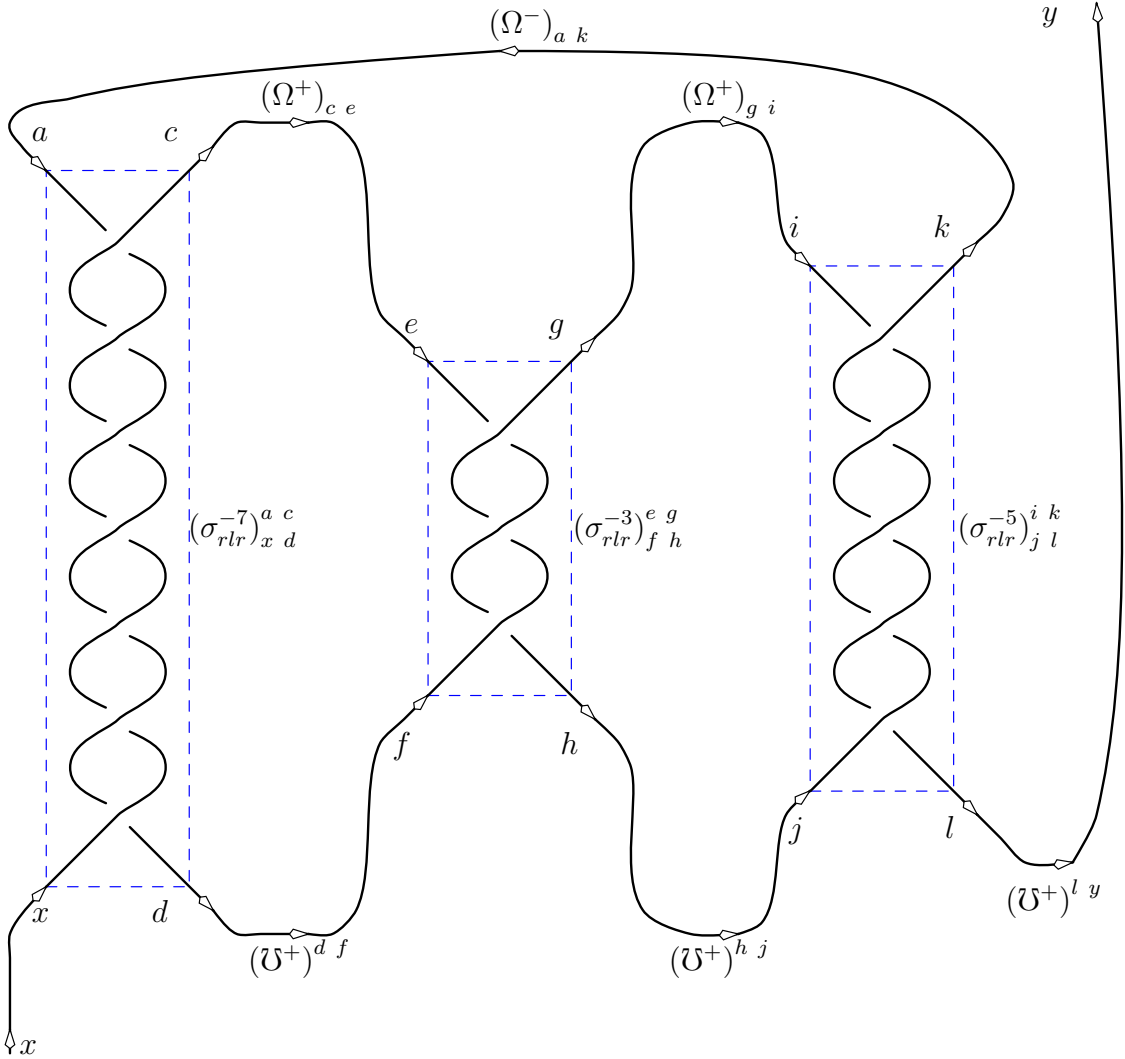}
    \caption{The (noninvertible) pretzel knots of
       Trotter, in tangle form. This illustration is of the smallest
       possible one, with $ p = 7 $, $ q = 3 $, $ r = 5 $.
    }
    \figlabel{TrotterPretzel}
  \end{center}
\end{figure}

The tensor associated with the pretzel is:
\be
  {T_{TP \( p, q, r \)}}^y_x
  & \defeq &
  {\( \sigma^{-p}_{rlr} \)}^{a~c}_{x~d}
  \cdot
  {\( \sigma^{-q}_{rlr} \)}^{e~g}_{f~h}
  \cdot
  {\( \sigma^{-r}_{rlr} \)}^{i~k}_{j~l}
  \cdot
  \\
  & & \qquad
  {\( \Omega^- \)}_{a~k}
  \cdot
  {\( \Omega^+ \)}_{c~e}
  \cdot
  {\( \Omega^+ \)}_{g~i}
  \cdot
  {\( \mho^+ \)}^{d~f}
  \cdot
  {\( \mho^+ \)}^{h~j}
  \cdot
  {\( \mho^+ \)}^{l~y}.
\ee
Experiments show that the Links--Gould invariant for this class of
noninvertible knots always displays the symmetry of
\eqref{Melanotaenia}, for all $ p, q, r \leqslant 67 $. This amounts to
$5456$ knots, the smallest being the $(3,5,7)$ pretzel, a knot of
$3+5+7=15$ crossings, and the largest being the $(63,65,67)$ pretzel, a
knot of $63+65+67=195$ crossings. Incidentally, we find that the
invariant demonstrates that all those pretzels are chiral.


\subsection{The Kinoshita--Terasaka Pair of Mutant Knots}
\seclabel{KTMutants}

The Kinoshita--Terasaka pair is an example of a pair of $ 11 $ crossing
mutant knots that are known to be distinct.  To be precise, more
commonly, the first of the pair is usually known as the
``Kinoshita--Terasaka Knot'', and the second has been called the
``Conway Knot''.  In the original source by Kinoshita and Terasaka
\cite[p~151]{KinoshitaTerasaka:57}, the knot involved is the one
labelled $ \kappa \( 2, 2 \) $ (reproduced in
\cite[p~53]{Livingston:93}).  They had constructed this knot as an
example of a nontrivial $ 11 $ crossing knot with Alexander polynomial
equal to $ 1 $.  The source used to draw our example is from
\cite[p~174]{Adams:94}; note that these diagrams have $ 12 $ crossings,
so they are not minimal.

A number of proofs of their distinctness are at hand:
\begin{itemize}
\item
  Adams \cite[p~106]{Adams:94} states that Francis Bonahon and Lawrence
  Siebenmann first showed this in 1981.
  Adams \cite[p~174]{Adams:94} goes on to state that
  in 1986 David Gabai \cite{Gabai:86} showed that
  their \emph{minimal genus Seifert surfaces} have different
  \emph{genera}%
  \footnote{%
    \textsl{
      If $ L $ is an oriented link in $ S^3 $
    }
    (i.e. the $ 3 $-sphere),
    \textsl{
      then a Seifert surface for $ L $ is an oriented surface $ R $
      embedded in $ S^3 $ such that $ \partial R = L $ and no component
      is closed.}
    \cite[p~677]{Gabai:86}.
    That is, a Seifert surface is a $ 2 $-manifold with boundary
    being the link in question; the genus of such a surface being a
    topological classifying label \cite[p~95-106]{Adams:94}. The
    original reference for the Seifert algorithm is contained in
    \cite{Seifert:34}.
  }.
\item
  More recently, Morton and Cromwell \cite{MortonCromwell:96} have
  constructed a Vassiliev invariant of \emph{type}%
  \footnote{%
    A Vassiliev invariant is defined \cite[p~229]{MortonCromwell:96} to
    be of \emph{type} $d$ if it is zero on any link diagram of
    $d+1$ nodes, and to be of \emph{degree} $d$ if it is of type $d$
    but not of type $d-1$.
  }
  $11$ which distinguishes them. This Vassiliev invariant is based on
  the HOMFLY polynomial for framed links, and the authors compare it
  with another invariant, itself coming from $ SU_q \( 3 \) $, which
  does \emph{not} distinguish them.

  More specifically, they show that the `$ SU_q \( N \) $ invariant'
  for the module with Young diagram
  $ \stackrel{\textstyle \Box \! \Box}{\Box \; \;} $ will distinguish
  at least \emph{some} mutant pairs (in particular the KT pair), for
  all $ N \geqslant 4 $, but will definitely \emph{not} distinguish any
  for $ N = 2, 3 $.
\end{itemize}

More generally, it is known that neither the HOMFLY nor the Kauffman
polynomial can distinguish \emph{any} pair of mutants
\cite[p~174]{Adams:94}. In fact Lickorish \cite{Lickorish:87} used
skein theoretical arguments to show this; and furthermore, Lickorish
and Lipson \cite{LickorishLipson:87} and, independently Przytycki
\cite{Przytycki:89} again used skein theoretical arguments to show that
two equally twisted \emph{$ 2 $-cables}
(definition in \cite[p~118]{Adams:94})
of a mutant pair would have the same HOMFLY polynomial. Perhaps the
strongest statement that can be made in this direction was provided in
1994 by Chmutov, Duzhin and Lando \cite{ChmutovDuzhinLando:94}, who
proved that \emph{all} Vassiliev invariants of type less than $ 9 $
will agree on \emph{any} pair of mutants.

The question of whether the Links--Gould invariant is able to
distinguish mutants is immediately answerable in the negative. Theorem
5 of \cite{MortonCromwell:96} states that if the modules occurring in
the decomposition of $ V \otimes V $ each have unit multiplicity, as
indeed \eqref{TPdecomp} shows in our case, then the invariant is unable
to detect mutations.  Whilst this was proved in
\cite{MortonCromwell:96} for the case of quantum algebras, the
extension to the case of quantum superalgebras is quite
straightforward. As an example, we have explicitly evaluated the
Links--Gould invariant for the aforementioned pair of mutants.

\pagebreak

We illustrate $ KT $, the first of Kinoshita--Terasaka pair, in
\figref{KT}, where the tensors $ KTA $ and $ KTB $ are
defined below, in Figures \ref{fig:KTA} and \ref{fig:KTC}. From $ KT $,
we may build the mutant $ KT' $ by replacing the component $ KTA $ with
$ KTA' $ (depicted in \figref{KTAprime}), which is formed by reflection
of $ KTA $ about a horizontal line.

\begin{figure}[htbp]
  \begin{center}
    \input{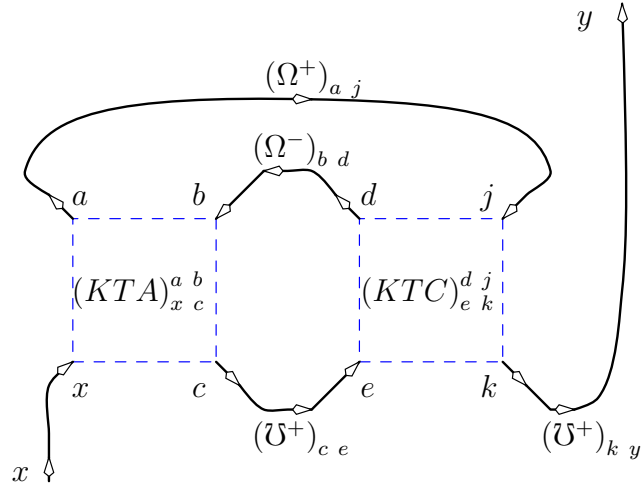}
    \caption{%
      $ KT $, the first of the Kinoshita--Terasaka pair of mutant knots,
      where the subdiagrams $ KTA $, $ KTA' $ and $ KTC $ are found in
      Figures \ref{fig:KTA}, \ref{fig:KTAprime} and \ref{fig:KTC}
      respectively.  (The mutant $ KT' $ of $ KT $ is obtained by
      exchanging $ KTA $ with $ KTA' $.)
    }
    \figlabel{KT}
  \end{center}
\end{figure}

The tensors associated with $ KT $ and $ KT' $ are:
\be
  {\( T_{KT} \)}^y_x
  & \defeq &
  {\( KTA \)}^{a~b}_{x~c}
  \cdot
  {\( KTC \)}^{d~j}_{e~k}
  \cdot
  {\( \Omega^- \)}_{b~d}
  \cdot
  {\( \mho^+ \)}^{c~e}
  \cdot
  {\( \Omega^+ \)}_{a~j}
  \cdot
  {\( \mho^+ \)}^{k~y},
  \\
  {\( T_{KT'} \)}^y_x
  & \defeq &
  {\( KTA' \)}^{a~b}_{x~c}
  \cdot
  {\( KTC \)}^{d~j}_{e~k}
  \cdot
  {\( \Omega^- \)}_{b~d}
  \cdot
  {\( \mho^+ \)}^{c~e}
  \cdot
  {\( \Omega^+ \)}_{a~j}
  \cdot
  {\( \mho^+ \)}^{k~y},
\ee
where
\be
  {\( KTA \)}^{a~b}_{q~c}
  & \defeq &
  {\( \sigma \sigma_d \)}^{a~b}_{d~e}
  \cdot
  {\( \sigma^{-2} \)}^{d~f}_{q~g}
  \cdot
  {\( \sigma^{-1}_d \)}^{h~e}_{i~c}
  \cdot
  {\( \Omega^+ \)}_{f~h}
  \cdot
  {\( \mho^- \)}^{g~i}
  \\
  {\( KTA' \)}^{a~b}_{q~c}
  & \defeq &
  {\( \sigma^{-2} \)}^{a~f}_{d~g}
  \cdot
  {\( \sigma^{-1}_d \)}^{h~b}_{i~e}
  \cdot
  {\( \sigma \sigma_d \)}^{d~e}_{q~c}
  \cdot
  {\( \Omega^+ \)}_{f~h}
  \cdot
  {\( \mho^- \)}^{g~i}
  \\
  {\( KTC \)}^{d~j}_{e~k}
  & \defeq &
  {\( KTB \)}^{d~f}_{e~g}
  \cdot
  {\( \sigma^{-1}_l \sigma^{-1}_r \)}^{h~j}_{i~k}
  \cdot
  {\( \Omega^- \)}_{f~h}
  \cdot
  {\( \mho^+ \)}^{g~i}
  \\
  {\( KTB \)}^{d~f}_{e~g}
  & \defeq &
  \sigma^{d~b}_{a~c}
  \cdot
  {( \sigma_d^2 )}^{l~f}_{m~n}
  \cdot
  {\( \sigma^{-1} \sigma^{-1}_d \)}^{a~n}_{e~g}
  \cdot
  {\( \Omega^+ \)}_{b~l}
  \cdot
  {\( \mho^- \)}^{c~m}.
\ee

\begin{figure}[htbp]
  \begin{center}
    \input{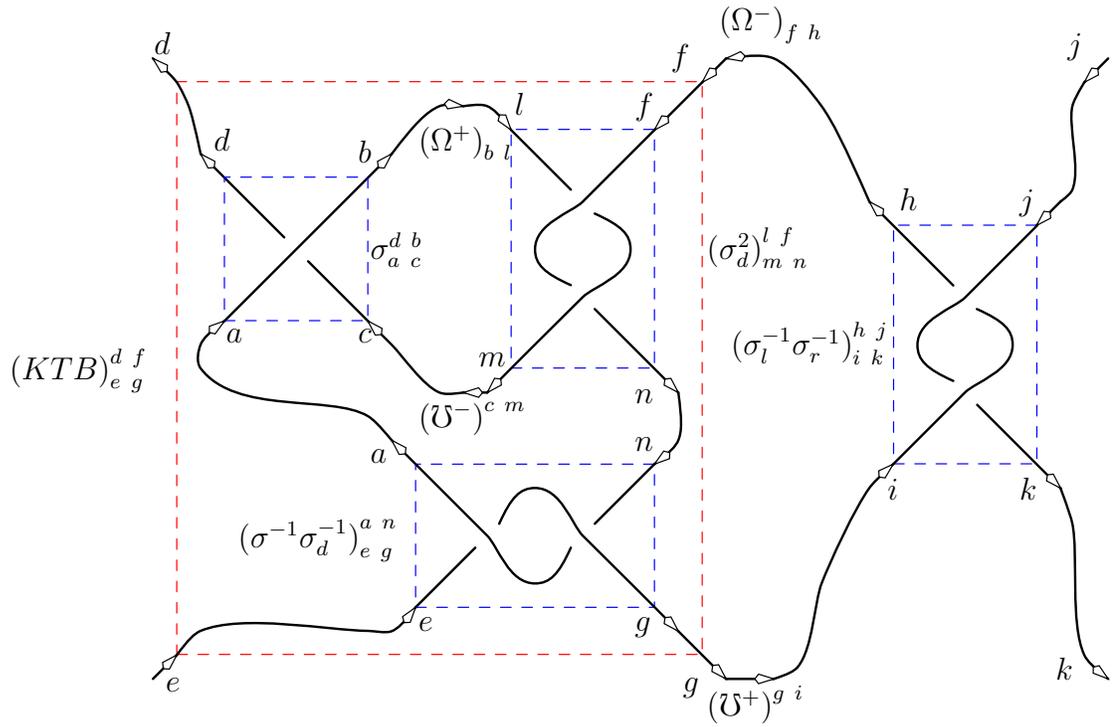}
    \caption{%
      The component $ KTC $ of the Kinoshita--Terasaka pair
      of mutant knots $ KT $ and $ KT' $.
    }
    \figlabel{KTC}
  \end{center}
\end{figure}

\begin{figure}[htbp]
  \begin{center}
    \input{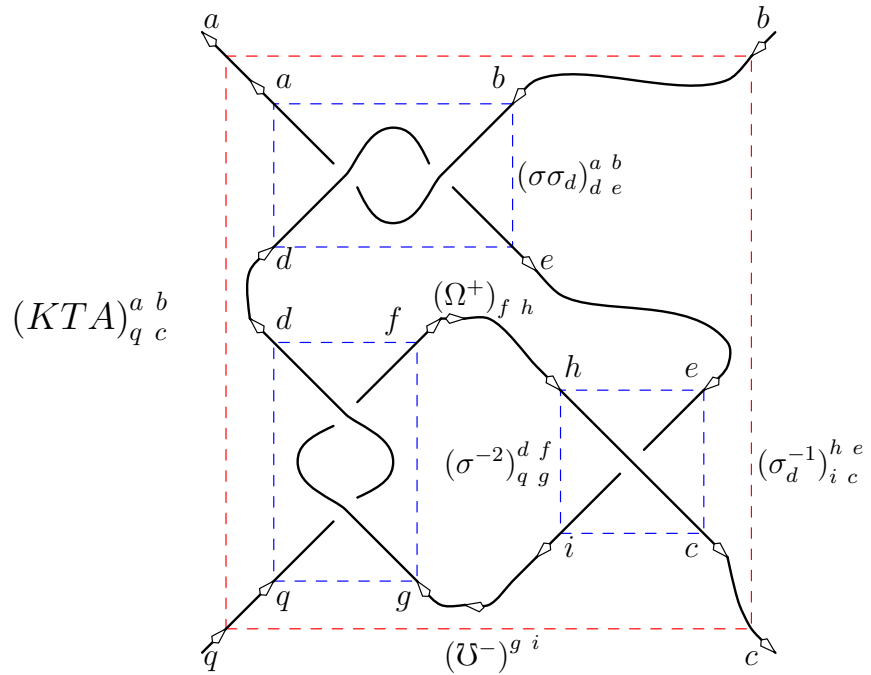}
    \caption{%
      The component $ KTA $ of $ KT $, the first of the
      Kinoshita--Terasaka pair.
    }
    \figlabel{KTA}
  \end{center}
\end{figure}

\begin{figure}[htbp]
  \begin{center}
    \input{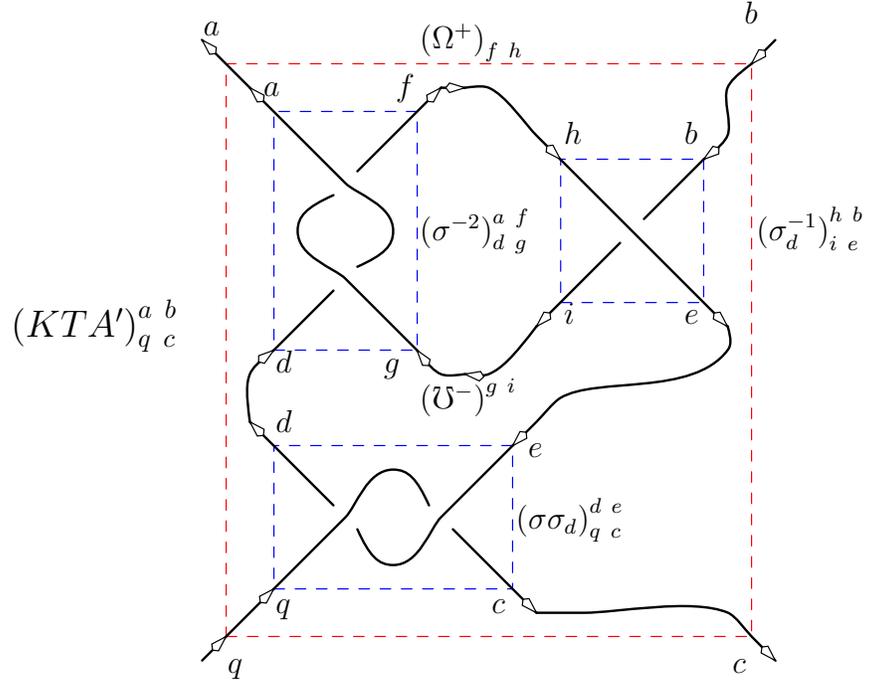}
    \caption{%
      The component $ KTA' $ of $ KT' $, the second of the
      Kinoshita--Terasaka pair.
    }
    \figlabel{KTAprime}
  \end{center}
\end{figure}

\pagebreak


\subsection{Links--Gould Polynomials of the K--T Mutants}

We find the Links--Gould polynomials of both mutants to be:

\begin{eqnarray*}
  & &
  \! \! \! \! \! \! \! \! \! \! \! \! \! \! \! \!
  \! \! \! \! \! \! \! \!
  {LG}_{KT} \( q, p \)
  =
  \\
  & &
  \! \! \! \! \! \! \! \! \! \! \! \! \! \! \! \!
  \begin{array}{l}
    - 23
    -    p^{-6} q^{-8}
    -    p^{-6} q^{-6}
    +  2 p^{-6} q^{-4}
    +    p^{-6} q^{-2}
    -    p^{-6}
    +    p^{-4} q^{-8}
    +  6 p^{-4} q^{-6}
    -  3 p^{-4} q^{-4}
    \\
    -  9 p^{-4} q^{-2}
    +  2 p^{-4}
    +  3 p^{-4} q^2
    -  7 p^{-2} q^{-6}
    -  7 p^{-2} q^{-4}
    + 18 p^{-2} q^{-2}
    +  9 p^{-2}
    - 11 p^{-2} q^2
    \\
    -  2 p^{-2} q^4
    +  2        q^{-6}
    + 14        q^{-4}
    -  8        q^{-2}
    +  6        q^2
    + 10        q^4
    -  7 p^2    q^{-4}
    -  7 p^2    q^{-2}
    + 18 p^2
    +  9 p^2    q^2
    \\
    - 11 p^2    q^4
    -  2 p^2    q^6
    +    p^4    q^{-4}
    +  6 p^4    q^{-2}
    -  3 p^4
    -  9 p^4    q^2
    +  2 p^4    q^4
    +  3 p^4    q^6
    -    p^6    q^{-2}
    -    p^6
    \\
    +  2 p^6    q^2
    +    p^6    q^4
    -    p^6    q^6,
  \end{array}
\end{eqnarray*}
hence the Links--Gould link invariant \emph{does not} distinguish
between these mutants.

As predicted by the theorem of \cite{MortonCromwell:96}, the tensors
$KTA$ and $KTA'$ are in fact \emph{identical}, which explains why the
pair of mutants yield the same invariant.

\pagebreak


\section*{Acknowledgements}

Louis Kauffman  thanks the National Science Foundation  for  support of
this research under grant number  DMS-9205277 and the NSA for partial
support under grant number MSPF-96G-179.

Jon Links is supported by an Australian Postdoctoral Research
Fellowship and a University of Queensland New Staff Research Grant.
The authors wish to thank Mark Gould for his enthusiastic support of
this research.


\def\JMP{Journal of Mathematical Physics}
\def\JPA{Journal of Physics A. Mathematical and General}
\def\LMP{Letters in Mathematical Physics}
\def\NPB{Nuclear Physics B}
\def\RMP{Reports on Mathematical Physics}
\def\CMP{Communications in Mathematical Physics}
\def\JKTR{Journal of Knot Theory and its Ramifications}



\bibliographystyle{plain}
\bibliography{lgl}

\begin{thebibliography}{10}

\bibitem{Adams:94}
Colin~C Adams.
\newblock {\em {The Knot Book: an Elementary Introduction to Mathematical
  Theory of Knots}}.
\newblock Freeman, New York, 1994.

\bibitem{Alexander:23}
James~W Alexander.
\newblock Topological invariants of knots and links.
\newblock {\em Transactions of the American Mathematical Society}, 20:275--306,
  1923.

\bibitem{AlexanderBriggs:26}
James~W Alexander and G~B Briggs.
\newblock On types of knotted curves.
\newblock {\em Annals of Mathematics}, 28:562--586, 1926-1927.

\bibitem{BarNatan:95}
Dror {Bar-Natan}.
\newblock On the {V}assiliev knot invariants.
\newblock {\em Topology}, 34(2):423--472, 1995.

\bibitem{ChmutovDuzhinLando:94}
S~V Chmutov, S~V Duzhin, and S~K Lando.
\newblock {V}assiliev knot invariants. {I}. introduction. singularities and
  bifurcations.
\newblock {\em Advances in Soviet Mathematics}, 21:117--126, 1994.

\bibitem{Conway:70}
John~H Conway.
\newblock An enumeration of knots and links, and some of their algebraic
  properties.
\newblock In Leech \cite{Leech:70}, pages 329--358.
\newblock Proceedings of a Conference in Oxford, 1967.

\bibitem{CrowellFox:77}
Richard~H Crowell and Ralph~H Fox.
\newblock {\em Introduction to Knot Theory}.
\newblock Number~57 in Graduate Texts in Mathematics. Springer-Verlag, Berlin,
  New York, {S}pringer edition, 1977.

\bibitem{DollHoste:91}
Helmut Doll and Jim Hoste.
\newblock A tabulation of oriented links.
\newblock {\em Mathematics of Computation}, 57(196):747--761, 1991.

\bibitem{Drinfeld:87}
V~G {Drinfel'd}.
\newblock Quantum groups.
\newblock In Andrew~M Gleason, editor, {\em Proceedings of the International
  Congress of Mathematicians 1986 (2 volumes)}, pages 798--820, Providence,
  Rhode Island, USA, 1987. American Mathematical Society.

\bibitem{FreydYetterHosteLickorishMilletOcneanu:85}
Peter Freyd, David~N Yetter, Jim Hoste, W~B~Raymond Lickorish, Kenneth~C
  Millet, and Adrian Ocneanu.
\newblock A new polynomial invariant of knots and links.
\newblock {\em Bulletin of the American Mathematical Society (New Series)},
  12(2):239--246, April 1985.
\newblock Included in \cite[pp~12-19]{Kohno:90}.

\bibitem{Gabai:86}
David Gabai.
\newblock Genera of the alternating links.
\newblock {\em Duke Mathematical Journal}, 53(3):677--681, 1986.

\bibitem{GouldLinksZhang:96b}
Mark~D Gould, Jon~R Links, and Yao-Zhong Zhang.
\newblock Type-{I} quantum superalgebras, {$ q $}-supertrace and two-variable
  link polynomials.
\newblock {\em \JMP}, 37:987--1003, 1996.

\bibitem{Jimbo:85}
Michio Jimbo.
\newblock A {$q$}-difference analogue of {$U(\mathfrak{g})$} and the
  {Y}ang-{B}axter equation.
\newblock {\em Letters in Mathematical Physics}, 10(1):63--69, 1985.

\bibitem{Jones:85}
Vaughan F~R Jones.
\newblock A polynomial invariant for knots via {Von Neumann} algebras.
\newblock {\em Bulletin of the American Mathematical Society}, 12(1):103--111,
  January 1985.
\newblock Included in \cite[pp~3-11]{Kohno:90}.

\bibitem{Jones:87}
Vaughan F~R Jones.
\newblock {H}ecke algebra representations of braid groups and link polynomials.
\newblock {\em Annals of Mathematicals}, 126:335--388, 1987.
\newblock Included in \cite[pp~20-73]{Kohno:90}.

\bibitem{Kauffman:87a}
Louis~H Kauffman.
\newblock {\em {On Knots}}.
\newblock Number 115 in Annals of Mathematics Studies. Princeton University
  Press, Princeton, NJ, 1987.

\bibitem{Kauffman:88}
Louis~H Kauffman.
\newblock New invariants in the theory of knots.
\newblock {\em The American Mathematical Monthly}, 95(3):195--242, March 1988.

\bibitem{Kauffman:93}
Louis~H Kauffman.
\newblock {\em {Knots and Physics}}.
\newblock World Scientific, Singapore, 2nd edition, 1993.

\bibitem{Kauffman:97a}
Louis~H Kauffman.
\newblock Knots and diagrams.
\newblock In Shin'ichi Suzuki, editor, {\em Lectures at {K}nots96}, pages
  123--194. World Scientific, 1997.

\bibitem{KhoroshkinTolstoy:91}
S~M Khoroshkin and V~N Tolstoy.
\newblock Universal {$ R $}-matrix for quantized (super)algebras.
\newblock {\em Communications in Mathematical Physics}, 141(3):599--617, 1991.

\bibitem{KinoshitaTerasaka:57}
Shin'ichi Kinoshita and Hidetaka Terasaka.
\newblock On unions of knots.
\newblock {\em Osaka Mathematical Journal}, 9:131--153, 1957.

\bibitem{Kohno:90}
Toshitake Kohno, editor.
\newblock {\em New Developments in the Theory of Knots}, volume~11 of {\em
  Advanced Series in Mathematical Physics}.
\newblock World Scientific, Singapore, 1990.

\bibitem{Leech:70}
John Leech, editor.
\newblock {\em Computational Problems in Abstract Algebra}, Oxford, UK, 1970.
  Pergamon Press.
\newblock Proceedings of a Conference in Oxford, 1967.

\bibitem{Lickorish:87}
W~B~Raymond Lickorish.
\newblock Linear skein theory and link polynomials.
\newblock {\em Topology and its Applications}, 27(3):265--274, 1987.

\bibitem{LickorishLipson:87}
W~B~Raymond Lickorish and Andrew~S Lipson.
\newblock Polynomials of $2$-cable-like links.
\newblock {\em Proceedings of the American Mathematical Society},
  100(2):355--361, 1987.

\bibitem{LinksGould:92b}
Jon~R Links and Mark~D Gould.
\newblock Two variable link polynomials from quantum supergroups.
\newblock {\em \LMP}, 26(3):187--198, November 1992.

\bibitem{LinksGouldZhang:93}
Jon~R Links, Mark~D Gould, and Rui~Bin Zhang.
\newblock Quantum supergroups, link polynomials and representation of the braid
  generator.
\newblock {\em Reviews in Mathematical Physics}, 5(2):345--361, 1993.

\bibitem{Livingston:93}
Charles Livingston.
\newblock {\em Knot Theory}, volume~24 of {\em The {C}arus Mathematical
  Monographs}.
\newblock The Mathematical Association of America, Washington, DC, 1993.

\bibitem{MortonCromwell:96}
Hugh~R Morton and Peter~R Cromwell.
\newblock Distinguishing mutants by knot polynomials.
\newblock {\em \JKTR}, 5(2):225--238, 1996.

\bibitem{PrasolovSossinsky:96}
Viktor~V Prasolov and A~B Sossinsky.
\newblock {\em Knots, Links, Braids and $ 3 $-Manifolds: An Introduction to the
  New Invariants in Low-Dimensional Topology}, volume 154 of {\em Translations
  of Mathematical Monographs}.
\newblock American Mathematical Society, Providence, Rhode Island, USA, 1996.
\newblock Translated from the Russian manuscript by Sossinsky.

\bibitem{Przytycki:89}
{J\'{o}zef}~H Przytycki.
\newblock Equivalence of cables of mutants of knots.
\newblock {\em Canadian Journal of Mathematics}, 41(2):250--273, 1989.

\bibitem{PrzytyckiTraczyk:87}
{J\'{o}zef}~H Przytycki and {Pawe\l} Traczyk.
\newblock Invariants of links of {C}onway type.
\newblock {\em Kobe Journal of Mathematics}, 4(2):115--139, 1987.

\bibitem{Reidemeister:48}
Kurt Reidemeister.
\newblock {\em Knotentheorie}.
\newblock Ergebnisse der Mathematik und ihrer Grenzgebiete; Bd 1 Hft 1.
  Chelsea, New York, 1948.
\newblock First published in German in 1932.

\bibitem{Reshetikhin:87}
Nikolai~Yu Reshetikhin.
\newblock Quantized universal enveloping algebras, the {Y}ang-{B}axter equation
  and invariants of links {I}, {II}.
\newblock L.O.M.I. preprints, numbers E-4-87 and E-17-87, never published,
  1987.

\bibitem{ReshetikhinTuraev:90}
Nikolai~Yu Reshetikhin and Vladimir~G Turaev.
\newblock Ribbon graphs and their invariants derived from quantum groups.
\newblock {\em \CMP}, 127:1--26, 1990.

\bibitem{Rolfsen:76}
Dale Rolfsen.
\newblock {\em Knots and Links}, volume~7 of {\em Mathematics Lecture Series}.
\newblock Publish or Perish, Inc, Wilmington, Delaware, USA, 1976.

\bibitem{Seifert:34}
Herbert Seifert.
\newblock Uber das {G}eschlecht von {K}noten.
\newblock {\em Mathematische Annalen}, 110:571--592, August 1934.
\newblock In German.

\bibitem{Trotter:64}
Hale~F Trotter.
\newblock Non-invertible knots exist.
\newblock {\em Topology}, 2:275--280, 1964.

\bibitem{Turaev:88}
Vladimir~G Turaev.
\newblock The {Y}ang-{B}axter equation and invariants of links.
\newblock {\em Inventiones Mathematicae}, 92(3):527--553, 1988.

\bibitem{Zhang:95}
Rui~Bin Zhang.
\newblock Quantum supergroups and topological invariants of three-manifolds.
\newblock {\em Review of Mathematical Physics}, 7(5):809--831, 1995.

\bibitem{ZhangGouldBracken:91b}
Rui~Bin Zhang, Mark~D Gould, and Anthony~J Bracken.
\newblock Quantum group invariants and link polynomials.
\newblock {\em \CMP}, 137(1):13--27, 1991.

\end{thebibliography}


\end{document}